\documentclass[a4paper,12pt]{amsart}
\usepackage{graphicx} 
\usepackage{inputenc}
\usepackage{tikz}
\usetikzlibrary{cd}
\usetikzlibrary{arrows.meta}
\usepackage{CJK}
\usepackage[colorlinks=true,linkcolor=blue,citecolor=blue,urlcolor=blue]{hyperref}
\usepackage{geometry,url,color}
\usepackage{amsfonts,amssymb,amsmath,fancyhdr,mathrsfs}
\usepackage{amsthm,array,ragged2e,graphicx,enumitem}
\usepackage{booktabs,float,multirow,diagbox,longtable}
\usepackage{xcolor}
\usepackage{capt-of}
\usepackage[numbers,square]{natbib}

\definecolor{mycolor}{HTML}{C7EDCC}    

\allowdisplaybreaks
\justifying
\numberwithin{equation}{section}
\newtheorem{theorem}{Theorem}[section]
\newtheorem{definition}[theorem]{Definition}
\newtheorem{lemma}[theorem]{Lemma}
\newtheorem{proposition}[theorem]{Proposition}

\newtheorem{corollary}[theorem]{Corollary}
\newtheorem{remark}[theorem]{Remark}

\newcommand{\Aut}{\mathrm{Aut}}
\newcommand{\Hom}{\mathrm{Hom}}

\newcommand{\id}{\mathrm{id}}
\newcommand{\re}{\mathrm{Re}}
\newcommand{\Ker}{\mathrm{Ker}}
\newcommand{\im}{\mathrm{Im}}

\newcommand{\diff}{\mathrm{d}}
\newcommand{\mmod}{\mathrm{mod}}

\newcommand{\rep}{\mathrm{rep}}

\newcommand{\Stab}{\mathrm{Stab}}

\newcommand{\linp}[1]{\left\langle#1\right\rangle}
\newcommand{\udv}{\underline{\dim}}

\newcommand{\supp}{\mathrm{supp}}

\newcommand{\sgn}{\mathrm{sgn}}

\newcommand{\Ss}{\mathrm{SS}}
\newcommand{\codim}{\mathrm{codim}}

\newcommand{\add}{\mathrm{add}}

\newcommand{\Iso}{\mathrm{Iso}}
\newcommand{\cl}{\mathrm{cl}}
\newcommand{\1}{\mathbf{1}}
\newcommand{\fin}{\mathrm{fin}}
\newcommand{\pr}{\mathrm{pr}}
\newcommand{\Arg}{\mathrm{Arg}}
\newcommand{\ev}{\mathrm{ev}}
\newcommand{\RHS}{\mathrm{RHS}}
\newcommand{\LHS}{\mathrm{LHS}}
\newcommand{\Lie}{\mathrm{Lie}}
\newcommand{\fl}{\mathrm{fl}}

\title[Numerical Joyce Transport and Derived Scattering Groupoid]{
Finite Numerical Joyce Transport and Derived Hall Scattering Groupoid
}
\author[Wenyu Gao and Fan Xu]{
Wenyu Gao* and Fan Xu
}
\address{Department of Mathematical Sciences\\
Tsinghua University\\
Beijing 100084, P.~R.~China} \email{gwy23@mails.tsinghua.edu.cn (W. Gao)}

\address{Department of Mathematical Sciences\\
Tsinghua University\\
Beijing 100084, P.~R.~China} \email{fanxu@mail.tsinghua.edu.cn (F. Xu)}
\keywords{Derived Hall algebra, Triangulated category, Stability condition, Scattering diagram, Joyce transport, Scattering groupoid}
	\thanks{$*$~Corresponding author.}

\begin{document}
\begin{abstract}
Let $\mathcal{D}$ be a triangulated category over a finite field $k$ satisftying certain finiteness conditions. Let $\Sigma$ be a connected component of the space of numerical stability conditions on $\mathcal{D}$ with one element $\sigma\in\Sigma$ satisfying the support property and the finite counting condition. We prove that each of the Joyce elements with respect to $\sigma_1\in\Sigma$ is a finite rational Lie polynomial in those of any other $\sigma_0\in\Sigma$, so the generated Lie subalgebra $\mathfrak{L}_\sigma,\sigma\in\Sigma$ in the uncompleted derived Hall algebra $\mathcal{H}(\mathcal{D})$ is constant on $\Sigma$. Our proof combines completions on strict sectors with a uniform mass bound for successive local substitutions. We use the method of Hall scattering diagrams on phase categories to provide a realization of the transform of Joyce elements. For the finite-length case, we also give an alternative proof using finite-length hearts instead of phase categories. As an application, we explicit the generators of the Lie subalgebra $\mathfrak{L}$ in the Dynkin case.
\end{abstract}
\maketitle
\section{Introduction}
    In \cite{Rei2003}, Reineke introduced the notion of Harder-Narasimhan filtration in representation categories of quivers, as an analog of a concept with the same name in the theory of moduli of vector bundles on curves. He established the relationship between stability functions and the Ringel-Hall algebras of representation categories. 

    Reineke's result was rewritten and summarized by Bridgeland and Toledano Laredo in \cite{BT2012} as the right part of the following diagram:
    
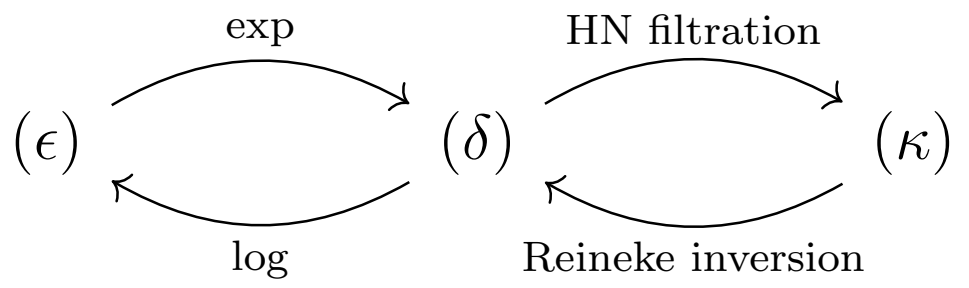
\begin{equation}\label{BT diagram}
    \begin{tikzcd}
(\epsilon) \arrow[rr, "\exp", bend left] &  & (\delta) \arrow[ll, "\log", bend left] \arrow[rr, "\text{HN filtration}", bend left] &  & (\kappa) \arrow[ll, "\text{Reineke inversion}", bend left]
\end{tikzcd}
\end{equation}
Roughly speaking, $(\kappa)$ and $(\delta)$ are characteristic functions, supported on all objects and semistable objects of a certain dimension vector, respectively.

In the left part of the diagram, the symbol $(\epsilon)$ denotes the elements defined by Joyce in \cite{Joy2007}, which will be called the Joyce elements in the present paper. These are not characteristic functions in general, but they behave well when the stability function changes. As shown by Joyce in \cite{Joy2008}, when the stability function changes under certain conditions, each Joyce element with respect to the new stability function can be expressed as a finite Lie polynomial in Joyce elements with respect to the old stability function. 

Let $\mathcal{A}$ be a finite length abelian category with finitely many isomorphism classes of simple objects. Reineke proved an identity \cite[Proposition 4.12]{Rei2003}, which is written as
\begin{equation}\label{Reineke wall-crossing}
    \1_\mathcal{A}=\prod^\leftarrow_{\phi\in(0,1]}\Ss^Z_\phi
\end{equation}
in the notation of Bridgeland and Toledano Laredo. Although it does not appear in the diagram, it plays a major role in the theory. It characterizes the existence and uniqueness of the Harder–Narasimhan filtration for a certain object in $\mathcal{A}$, while implying that the descending product on the right-hand side of the identity is an invariant under the change of stability functions.
    
Reineke's work was generalized by Joyce in \cite{Joy2007} and \cite{Joy2008}, where Joyce raised the question \cite[Problem 7.1]{Joy2008} of introducing a parallel theory on triangulated categories, based on the notion of stability conditions introduced by Bridgeland \cite{Bri2007}, which can be defined on bounded derived categories of abelian categories. Regarding Hall algebras, T\"{o}en \cite{Toe2006} constructed an associative algebra for the derived category of the category of perfect modules over a locally finite dg category, called the derived Hall algebra. This construction was later extended by Xiao and Xu \cite{XX2008} to triangulated categories satisfying certain conditions. They also provided a direct proof of the associativity of the derived Hall algebras.

With the preparation above, we are now ready to generalize Reineke's work to triangulated categories satisfying both the conditions of Bridgeland and those of Xiao and Xu. The categories satisfying all the conditions do exist. 

In such a category, Wang and Chen \cite{WC2023} investigated the HN filtrations  using the language of derived Hall algebras. Their result is similar in form to that of Reineke. The original version of the present paper \cite{GX2025} followed their approach and tried to recover the entirety of Reineke's theory and the diagram (\ref{BT diagram}). 

In the present paper, we will not try to recover the results on the whole triangulated category. Instead, we will work on each small segment under a stability condition and focus on the behavior of Joyce elements when the stability condition changes. The result of this paper might be a partial answer to the question \cite[Problem 7.1]{Joy2008} raised by Joyce.

\subsection{Connections with existing works}
Kontsevich and Soibelman prove finite Lie transport for continuous families of stability data on a graded Lie algebra \cite[Section 2.3, Corollary 1]{KS2008}. Their symmetric stability data record real phases and shifts. Over a finite field, they also impose finiteness at fixed class and phase, and identify the images of their sector factors in To\"en's derived Hall algebra \cite[Section 2.2, Definition 2; Section 6.1, Proposition 12]{KS2008}. In this paper, we work directly in the derived Hall algebra. For each Joyce element, we use a uniform mass bound and the finite counting condition to show that its transport involves only finitely many objects. We also describe each local step using a Hall scattering diagram on a phase category.

Mozgovoy \cite[Theorem 5.3]{Moz2022} proves an HN factorization in the completed motivic Hall algebra of a sector, using classes of stacks of semistable objects. If the stability condition has global dimension at most $2$, he also obtains a wall-crossing formula in a quantum torus \cite[Sections 5.4--5.6]{Moz2022}. In this paper, we work with functions on isomorphism classes over $\mathbb{F}_q$. We give a criterion for continuous Hall multiplication and bound the masses of the objects in each substitution. The resulting transport formula holds in the uncompleted derived Hall algebra, without an integration map.

\subsection{The main results of the paper}
Let $\mathcal{D}$ be a triangulated category over $\mathbb{F}_q$ satisfying the finiteness assumptions in Section 2. Theorem \ref{obstruction} shows that Hall multiplication has no jointly continuous extension to
\[\mathscr{F}(\Iso(\mathcal{D}))=\prod_{X\in\Iso(\mathcal{D})}\mathbb{Q}u_X,\]
where each factor has the discrete topology. Thus the global completion used in the earlier version \cite{GX2025} must be replaced. Nevertheless, the descending product
\[\1_{\mathcal{D}}=\prod_{\phi\in\mathbb{R}}^{\leftarrow}\1_{\mathcal{P}_\sigma(\phi)}\]
is valid coefficientwise by Theorem \ref{descending product}. Each coefficient in it is determined by the finite HN filtration of the corresponding object. It is a generalization of (\ref{Reineke wall-crossing}).

We work with strict sectors $\mathcal{P}_\sigma(I)$, where $I$ is contained in an interval of length less than $1$. Fix a surjective class map $\cl:K_0(\mathcal{D})\rightarrow\Gamma$ to a finite-rank lattice. We assume that $\sigma$ is numerical with respect to $\cl$ and satisfies the support property and the finite counting condition. The latter means that there are only finitely many isomorphism classes of semistable objects of each fixed phase and numerical class. These assumptions give a completed Hall algebra on each strict sector. In this algebra, we recover the HN identity, Reineke inversion, and the exponential and logarithmic formulas (\ref{HN identity in sector})-(\ref{Joyce inversion}). The resulting Joyce elements $\epsilon^\sigma_{\phi,\gamma}$ belong to $\mathcal{H}(\mathcal{D})$ and are independent of the sector used to define them.

The main theorem of this paper is Theorem \ref{Joyce transport}. Let $\Sigma$ be a connected component of $\Stab_\cl(\mathcal{D})$ containing a stability condition with the two properties above. For $\sigma_0,\sigma_1\in\Sigma$, every $\epsilon^{\sigma_1}_{\phi,\gamma}$ is a finite Lie polynomial with rational coefficients in the elements $\epsilon^{\sigma_0}_{\psi,\beta}$. Each Lie monomial has total numerical class $\gamma$. In particular,
\[\mathfrak{L}_{\sigma_0}=\mathfrak{L}_{\sigma_1}\subseteq\mathcal{H}(\mathcal{D}),\]
where $\mathfrak{L}_\sigma$ is the Lie subalgebra generated by all Joyce elements at $\sigma$. The proof combines local Joyce transformations with the ordered factorization of graded Lie algebras in \cite{KS2008}. A uniform mass bound along a path ensures that only finitely many substitutions are needed for each target element.

Following Bridgeland's construction \cite{Bri2017}, we realize the local transformations by consistent Hall scattering diagrams on the phase categories $\mathcal{P}_\sigma(\phi)$ at finite height cutoffs. The homogeneous components of their wall logarithms recover the Joyce elements for nearby stability conditions. We record successive substitutions by transport trees and forests. After identifying forests whose Lie polynomials have the same values in $\mathcal{H}(\mathcal{D})$, we obtain a transport category. Its isomorphisms form the scattering transport groupoid.

Finally, Theorem \ref{finite length Joyce transport} gives an alternative proof along paths of rotated algebraic stability conditions, using finite-length hearts. For a Dynkin quiver $Q$, every Joyce element in $\mathcal{H}(D^b(\rep_{\mathbb{F}_q}Q))$ is a finite rational Lie polynomial in
\[u_{S_i^{\oplus r}[b]},i\in Q_0,r\geqslant1,b\in\mathbb{Z},\]
where $S_i$ is the simple representation at the vertex $i$.

\subsection{Organization of the paper}
The paper is organized as follows. 

In Section 2, we recall the notion of derived Hall algebras, Bridgeland stability conditions and ordered factorization for graded Lie algebras. We also prove some basic results as preliminaries of Sections 4 and 6.

In Section 3, we rule out a jointly continuous extension of Hall multiplication to $\mathscr{F}(\Iso(\mathcal{D}))$ and use HN truncations to show that the descending product over phases equals $\1_{\mathcal{D}}$ coefficientwise.

In Section 4, we prove finite numerical Joyce transport by combining local factorizations in strict sectors. 

In Section 5, we realize local Joyce transport by consistent Hall scattering diagrams on phase categories and define the scattering transport groupoid using admissible transport forests between finite frames.

In Section 6, we use finite-length hearts to prove Joyce transport along paths in the space of rotated algebraic stability conditions and provide a precise calculation for Dynkin quivers as an application of the results.

\subsection{Conventions and notation}
    In this paper, every subcategory is assumed to be full. For a category $\mathcal{C}$, we write $X\in\mathcal{C}$ to mean that $X$ is an object in $\mathcal{C}$, and $\mathcal{B}\subseteq\mathcal{C}$ to mean that $\mathcal{B}$ is a subcategory of $\mathcal{C}$.

    For a category $\mathcal{C}$, we use $K^b(\mathcal{C})$ to denote the bounded homotopy category of $\mathcal{C}$. For an abelian category $\mathcal{A}$, we use $D^b(\mathcal{A})$ to denote the bounded derived category of $\mathcal{A}$.

    Let $\mathcal{S}$ be a subset of $\mathcal{C}$. We denote by $\add\mathcal{S}$ the additive hull of $\mathcal{S}$, namely the smallest additive subcategory of $\mathcal{C}$ containing all direct summands of finite direct sums of objects in $\mathcal{S}$. If $\mathcal{S}=\{X_1,\cdots,X_n\}$ is finite, we also write $\add(X_1,\cdots,X_n)=\add\mathcal{S}$.

    For a triangulated category $\mathcal{D}$, we use $[1]$ to denote the shift functor of $\mathcal{D}$. For a triangle
    \[X\rightarrow Y\rightarrow Z\rightarrow X[1],\]
    we may write it as
    \begin{center}
        \begin{tikzcd}
    X \arrow[rr] &                      & Y \arrow[ld] \\
                 & Z \arrow[lu, dashed] &             
    \end{tikzcd}
    \end{center}
    when necessary. For a heart of $\mathcal{D}$, we always mean that it is a heart of a bounded $t$-structure on $\mathcal{D}$.
\section{Derived Hall algebras and stability conditions}
\subsection{Derived Hall algebras}
Let $k=\mathbb{F}_q$ be a finite field with $q$ elements, and let $\mathcal{D}\neq0$ be an essentially small, idempotent-complete $k$-linear triangulated category satisfying the following finiteness conditions as in \cite{XX2008}:
\begin{itemize}
    \item $\mathcal{D}$ is Hom-finite;
    \item $\mathcal{D}$ is Krull-Schmidt, and
    \item $\mathcal{D}$ is left homologically finite, that is,
    \begin{equation}\label{eq:leftfinite}
 \sum_{i\geq0}\dim_k\Hom_{\mathcal{D}}(A[i],B)<\infty
 \quad\text{for all }A,B\in\mathcal{D}.
\end{equation}
\end{itemize}

Write $\Iso(\mathcal{D})$ for its set of isomorphism classes. We identify an object with its isomorphism class when this causes no ambiguity. For $A,B\in\mathcal{D}$, put
\[
 \{A,B\}=\prod_{i>0}|\Hom_{\mathcal{D}}(A[i],B)|^{(-1)^i},
 a_B=|\Aut_{\mathcal{D}}(B)|.
\]
The set of maps $X\rightarrow B$ whose cone is isomorphic to $A[1]$ is denoted by $\Hom(X,B)_{A[1]}$. Let
\[\mathcal{H}(\mathcal{D})=\bigoplus_{X\in\Iso(\mathcal{D})}\mathbb{Q}u_X,\]
and we define
\begin{equation}\label{eq:Hall}
 u_A*u_B=\sum_X F^X_{AB}u_X,\qquad
 F^X_{AB}=\frac{|\Hom(X,B)_{A[1]}|}{a_B}
              \frac{\{X,B\}}{\{B,B\}}.
\end{equation}
The rational number $F^X_{AB}$ is called the \emph{derived Hall number} with respect to the objects $A,B$ and $X$. 
By \cite[Theorem 3.6]{XX2008}, $(\mathcal{H}(\mathcal{D}),*)$ is an associative $\mathbb{Q}$-algebra with the identity $1=u_0$, called the \emph{derived Hall algebra} of $\mathcal{D}$.

For an element
\[f=\sum_{i=1}^k c_i u_{E_i}\in\mathcal{H}(\mathcal{D}),\]
we can view $f$ as a $\mathbb{Q}$-valued function on $\Iso(\mathcal{D})$:
\begin{align*}
    f:\Iso(\mathcal{D})&\rightarrow\mathbb{Q}\\
    E_i&\mapsto c_i,\\
    E_i\neq E&\mapsto 0.
\end{align*}
We will frequently use this point of view in this article. It is clear that any element in $\mathcal{H}(\mathcal{D})$, as a function, is finitely supported.

Let $\mathcal{C}$ be a subcategory of $\mathcal{D}$, with $\Iso(\mathcal{C})$ finite. We denote by $\1_{\mathcal{C}}$ the element in $\mathcal{H}(\mathcal{D})$ taking the value $1$ on objects in $\mathcal{C}$ and $0$ elsewhere. 

We write $K_0(\mathcal{D})$ for the \emph{Grothendieck group} of $\mathcal{D}$, whose element corresponding to an object $E\in\mathcal{D}$ is denoted by $\udv E$. $\mathcal{H}(\mathcal{D})$ has a natural grading with respect to $K_0(\mathcal{D})$:
\[\mathcal{H}(\mathcal{D})=\bigoplus_{\gamma\in K_0(\mathcal{D})}\mathcal{H}(\mathcal{D})_\gamma,\]
where $\mathcal{H}(\mathcal{D})_\gamma$ is the subspace of $\mathcal{H}(\mathcal{D})$ spanned by $u_E$ with $E\in\mathcal{D}$ and $\udv E=\gamma$.
\subsection{Bridgeland stability conditions}
\begin{definition}
    Let $\mathcal{D}$ be a triangulated category. A \emph{Bridgeland stability condition} on $\mathcal{D}$ is a pair $\sigma=(Z,\mathcal{P})$, where
    \begin{itemize}
        \item $Z:K_0(\mathcal{D})\rightarrow\mathbb{C}$ is a group homomorphism, called the \emph{central charge};
        \item $\mathcal{P}=\{\mathcal{P}(\phi)|\phi\in\mathbb{R}\}$ is a family of full additive extension-closed subcategories of $\mathcal{D}$,
    \end{itemize}
    satisfying the following conditions:
    \begin{enumerate}
        \item If $0\neq M\in\mathcal{P}(\phi)$, then $Z(M)=m_Me^{i\pi\phi}$ for some $m_M\in\mathbb{R}_+$.
        \item $\mathcal{P}(\phi+1)=\mathcal{P}(\phi)[1]$ for any $\phi\in\mathbb{R}$.
        \item For $A_i\in\mathcal{P}(\phi_i),i=1,2$ with $\phi_1>\phi_2$, we have $\Hom_\mathcal{D}(A_1,A_2)=0$.
        \item (HN filtration) For any $0\neq E\in\mathcal{D}$, there is a sequence of triangles
        \begin{center}
            \begin{tikzcd}[sep=tiny]
0=E_0 \arrow[rr] &                        & E_1 \arrow[rr] \arrow[ld] &                        & E_2 \arrow[r] \arrow[ld] & \cdots \arrow[r] & E_{n-1} \arrow[rr] &                        & E_n=E \arrow[ld] \\
                 & A_1 \arrow[lu, dashed] &                           & A_2 \arrow[lu, dashed] &                          &                  &                    & A_n \arrow[lu, dashed] &                 
\end{tikzcd}
        \end{center}
        such that each $A_i\neq0$ belongs to $\mathcal{P}(\phi_i)$ for some $\phi_i$, $1\leqslant i\leqslant n$, with $\phi_1>\cdots>\phi_n$. This sequence is called a \emph{Harder-Narasimhan filtration}, or an \emph{HN filtration} for short. It is unique up to isomorphism by \cite{Bri2007}. Thus we can define
        \[\phi^+(E)=\phi_1,\phi^-(E)=\phi_n,\Delta(E)=\phi_1-\phi_n\]
        and
        \[m(E):=\sum_{j=1}^n|Z(A_j)|.\]
    \end{enumerate}
\end{definition}

For an interval $I\subseteq\mathbb{R}$, let $\mathcal{P}(I)$ consist of zero and the objects whose HN factors have phases in $I$, that is,
\[\mathcal{P}(I)=\{0\neq E\in\mathcal{D}|\phi^{\pm}(E)\in I\}\cup\{0\}.\]
It is extension-closed. The subcategory $\mathcal{P}(0,1]$ is called the \emph{heart} of $\sigma$. It is an abelian category by \cite[Section 3]{Bri2007}. For any $\phi\in\mathbb{R}$, the subcategory $\mathcal{P}(\phi,\phi+1]$ is also abelian and called a \emph{rotated heart}. Bridgeland proved the following important result:
\begin{lemma}[{\cite[Proposition 5.3]{Bri2007}}]\label{heart}
    To give a stability condition on a triangulated category $\mathcal{D}$ is equivalent to giving a bounded t-structure on $\mathcal{D}$ and a stability function having the HN property on its heart.
\end{lemma}

By \cite[Lemma 5.2]{Bri2007}, the subcategories $\mathcal{P}(\phi),\phi\in\mathbb{R}$ are abelian. We call $\mathcal{P}(\phi)$ the \emph{phase category} with phase $\phi$. It will play a major role in our study. We say an object $0\neq E\in\mathcal{D}$ is \emph{$\sigma$-semistable} if $E$ belongs to some phase category $\mathcal{P}(\phi)$ for some phase $\phi\in\mathbb{R}$, and we say it is \emph{$\sigma$-stable} if it is simple in the abelian phase category $\mathcal{P}(\phi)$. For a $\sigma$-semistable object $E\in\mathcal{P}(\phi)$, we often write $\phi(E)=\phi$.

The following property is useful.
\begin{proposition}[Phase truncation]\label{phase truncation}
    Let $E\in\mathcal{D}$ and $c\in\mathbb{R}$. There is a unique (up to unique isomorphism of triangles inducing the identity on $E$) \emph{phase truncation triangle}
    \[E_{>c}\rightarrow E\rightarrow E_{\leqslant c}\rightarrow E_{>c}[1],\]
    where $E_{>c}\in\mathcal{P}(c,+\infty)$ and $E_{\leqslant c}\in\mathcal{P}(-\infty,c]$.
\end{proposition}
\begin{proof}
    If $E\in\mathcal{P}(-\infty,c]$, take $E_{\leqslant c}=E$ and $E_{>c}=0$. If $E\in\mathcal{P}(c,+\infty)$, take $E_{\leqslant c}=0$ and $E_{>c}=E$. Otherwise, $\phi^-(E)\leqslant c<\phi^+(E)$. Let
    \begin{center}
            \begin{tikzcd}[sep=tiny]
0=E_0 \arrow[rr] &                        & E_1 \arrow[rr] \arrow[ld] &                        & E_2 \arrow[r] \arrow[ld] & \cdots \arrow[r] & E_{n-1} \arrow[rr] &                        & E_n=E \arrow[ld] \\
                 & A_1 \arrow[lu, dashed] &                           & A_2 \arrow[lu, dashed] &                          &                  &                    & A_n \arrow[lu, dashed] &                 
\end{tikzcd}
        \end{center}
        be the HN filtration of $E$, with $\phi(A_1)>\phi(A_2)>\dots>\phi(A_{i-1})>c\geqslant\phi(A_i)>\dots>\phi(A_n)$, where $i\in\{2,\dots,n\}$ is uniquely determined. Take $E_{>c}=E_{i-1}$ and $E_{\leqslant c}$ the mapping cone of the composite morphism $E_{i-1}\rightarrow E_i\rightarrow\dots\rightarrow E$. It is extended from $A_i,\dots,A_n$ in order by applying the octahedral axiom several times. Thus we have completed the proof of existence.

        For the uniqueness part, let
        \[A\xrightarrow{i} E\xrightarrow{p} B\xrightarrow{d} A[1],A'\xrightarrow{i'} E\xrightarrow{p'} B'\xrightarrow{d'} A'[1]\]
        be two such phase truncations. Applying the cohomological functor $\Hom_{\mathcal{D}}(A,-)$ to the second triangle, we obtain the long exact sequence
        \[\dots\rightarrow\Hom(A,B'[-1])\rightarrow\Hom(A,A')\xrightarrow{i'\circ-}\Hom(A,E)\rightarrow\Hom(A,B')\rightarrow\cdots.\]
        By assumption, $B'\in\mathcal{P}(-\infty,c]$, so $B'[-1]\in\mathcal{P}(-\infty,c-1]$, but $A\in\mathcal{P}(c,+\infty)$, so
        \[\Hom(A,B'[-1])=\Hom(A,B')=0,\]
        and hence $i'\circ-$ is an isomorphism. So there is a unique morphism $\alpha:A\rightarrow A'$, such that $i'\circ\alpha=i$. Similarly, applying $\Hom_{\mathcal{D}}(A',-)$ to the first triangle gives that there is a unique morphism $\alpha':A'\rightarrow A$ such that $i\circ\alpha'=i'$. Putting the two composition equations together gives
        \[i\circ(\alpha'\circ\alpha)=i.\]
        By applying $\Hom_{\mathcal{D}}(A,-)$ to the first triangle, we deduce that $\alpha'\circ\alpha=\id_A$. Similarly, we have $\alpha\circ\alpha'=\id_{A'}$. Hence $\alpha$ and $\alpha'$ are isomorphisms between $A$ and $A'$, and they are unique respectively. By the axioms of triangulated categories, they can be extended to isomorphisms of triangles $(\alpha,\id_E,\beta)$ and $(\alpha',\id_E,\beta')$. By the same argument, we know $\beta$ and $\beta'$ are also uniquely determined. Thus we finish the proof of uniqueness.
\end{proof}
By the same argument, we obtain another truncation
\[E_{\geqslant c}\rightarrow E\rightarrow E_{<c}\rightarrow E_{\geqslant c}[1].\]
Hence, for any segment $I\subseteq\mathbb{R}$, by taking truncations for $E$ once or twice, we can find a unique object $E_I\in\mathcal{P}(I)$. We call $E_I$ the \emph{canonical HN block} of $E$ in $I$.

\begin{definition}[{\cite[Definition 5.7]{Bri2007}}]
    We say a stability condition $\sigma=(Z,\mathcal{P})$ on $\mathcal{D}$ is \emph{locally finite} if there exists $\eta>0$ such that $\mathcal{P}(\phi-\eta,\phi+\eta)$ is of finite length for any $\phi\in\mathbb{R}$.
\end{definition}
In the rest of this article, we will always assume that a stability condition is locally finite. As a consequence, the phase categories $\mathcal{P}(\phi),\phi\in\mathbb{R}$ are always of finite length.
This provides the foundation of our study in Sections 4-5.

Let $\Stab(\mathcal{D})$ be the set of locally finite stability conditions on $\mathcal{D}$. We introduce the generalized metric, that is, a metric that allows infinite distance, on $\Stab(\mathcal{D})$ defined by Bridgeland:
\begin{lemma}[{\cite[Proposition 8.1]{Bri2007}}]
    Let $\sigma,\tau\in\Stab(\mathcal{D})$. The function
    \[d(\sigma,\tau)=\sup_{0\neq E\in\mathcal{D}}\max\left\{|\phi_\sigma^-(E)-\phi_\tau^-(E)|,|\phi_\sigma^+(E)-\phi_\tau^+(E)|,\left|\log\dfrac{m_\sigma(E)}{m_\tau(E)}\right|\right\}\]
    defines a generalized metric on $\Stab(\mathcal{D})$.
\end{lemma}
We will use the resulting phase and mass bounds to compare nearby stability conditions.

The generalized metric $d$ induces a topology on $\Stab(\mathcal{D})$. Let $\Sigma$ be a connected component. Bridgeland defined a subspace $V(\Sigma)\subseteq\Hom_{\mathbb{Z}}(K_0(\mathcal{D}),\mathbb{C})$ and constructed a norm on it; see \cite[Section 6]{Bri2007}. Bridgeland then proved
\begin{proposition}[{\cite[Theorem 1.2]{Bri2007}}]\label{deformation}
    The map
    \begin{align*}
        \mathcal{Z}:\Sigma&\rightarrow V(\Sigma)\\
        \sigma=(Z,\mathcal{P})&\mapsto Z
    \end{align*}
    is a local homeomorphism. In particular, if $K_0(\mathcal{D})=\mathbb{Z}^n$ for some $n\in\mathbb{Z}_+$, then $\Stab(\mathcal{D})$ is a complex manifold.
\end{proposition}
If $\mathcal{D}=D^b(\mathcal{A})$ for some abelian category $\mathcal{A}$ of finite length with only finitely many isomorphism classes of simple objects, then $K_0(\mathcal{D})$ is a finitely generated free abelian group and hence $\Stab(\mathcal{D})$ is a complex manifold. However, in this article, we will not only study such categories. That motivates us to introduce the following concepts.

\subsection{Numerical stability conditions}
Let us fix a group epimorphism
\[\cl:K_0(\mathcal{D})\rightarrow\Gamma,\]
where $\Gamma$ is a finite rank lattice. For an object $E\in\mathcal{D}$, we always write $\cl E=\cl(\udv E)$ for short.
\begin{definition}
    A stability condition $\sigma=(Z,\mathcal{P})$ is called \emph{numerical} with respect to $\cl$ if $Z$ factors through $\cl$. The subspace of $\Stab(\mathcal{D})$ consisting of numerical stability conditions is denoted by $\Stab_\cl(\mathcal{D})$.
\end{definition}
The following is a natural consequence of Proposition \ref{deformation}.
\begin{corollary}[{\cite[Corollary 1.3]{Bri2007}}]
    The subspace $\Stab_\cl(\mathcal{D})$ is a complex manifold.
\end{corollary}
Fix a norm $\|\cdot\|$ on the complex vector space $\Gamma_{\mathbb{C}}=\Gamma\otimes_{\mathbb{Z}}\mathbb{C}$. We introduce the following concepts:
\begin{definition} Let $\sigma=(Z,\mathcal{P})\in\Stab_\cl(\mathcal{D})$.
    \begin{enumerate}
        \item $\sigma$ is said to have the \emph{support property} if there exists a $c>0$ such that
        \[|Z(\cl E)|\geqslant c\|\cl E\|\]
        for any $\sigma$-semistable object $E$.
        \item $\sigma$ is said to have the \emph{finite counting condition} if for any $\phi\in\mathbb{R}$ and $\gamma\in\Gamma$, the set
    \[\{E\in\Iso(\mathcal{D})|E\in\mathcal{P}(\phi),\cl E=\gamma\}\]
    is finite.
    \end{enumerate}
\end{definition}
We have the following properties.
\begin{proposition}\label{mass and phase bounds finite}
    Let $\sigma=(Z,\mathcal{P})\in\Stab_\cl(\mathcal{D})$ satisfy the support property and the finite counting condition. For any $R>0$ and any bounded segment $I\subseteq\mathbb{R}$, the set
    \[\{0\neq E\in\Iso(\mathcal{P}(I))|m(E)\leqslant R\}\]
    is finite.
\end{proposition}
\begin{proof}
    Let $A_1,\dots,A_r$ be the HN factors of $E$. Write $\gamma_j=\cl A_j\in\Gamma$. Let $c>0$ be the support constant and set
    \[\lambda=\min_{0\neq\gamma\in\Gamma}\|\gamma\|.\]
    By the discreteness of $\Gamma$, $\lambda>0$. The support property gives
    \[|Z(\gamma_j)|\geqslant c\|\gamma_j\|\geqslant c\lambda,j=1,\dots,r,\]
    while the upper bound of mass forces
    \[\sum_{j=1}^r|Z(\gamma_j)|=m(E)\leqslant R.\]
    Thus, we obtain
    \[r\leqslant\left[\dfrac{R}{c\lambda}\right],\|\gamma_j\|\leqslant\dfrac{R}{c},j=1,\dots,r.\]
    Since $\Gamma$ is a lattice, the bounds above give only finitely many choices for $\gamma_1,\dots,\gamma_r$. For each fixed $\gamma_j$, the phase of $A_j$ is determined modulo $2$ by $Z(\gamma_j)$. Since $E\in\mathcal{P}(I)$, this phase lies in $I$, so it has only finitely many possible values. The finite counting condition now gives only finitely many choices for each $A_j$ up to isomorphism. 
    
    Finally, for fixed $A,B\in\mathcal{D}$, the middle term of a triangle $A\rightarrow X\rightarrow B\rightarrow A[1]$ is determined up to isomorphism by a morphism $B\rightarrow A[1]$. There are only finitely many such morphisms because $\mathcal{D}$ is Hom-finite over $k=\mathbb{F}_q$. Applying this at each step of the HN filtration gives only finitely many isomorphism classes of $E$.
\end{proof}
\begin{corollary}\label{finite distance}
    Let $\sigma=(Z_\sigma,\mathcal{P}_\sigma)\in\Stab_\cl(\mathcal{D})$ satisfy the support property and the finite counting condition. If $\tau=(Z_\tau,\mathcal{P}_\tau)\in\Stab_\cl(\mathcal{D})$ with $d(\sigma,\tau)<+\infty$, then $\tau$ satisfies the support property and the finite counting condition.
\end{corollary}
\begin{proof}
    Let us first prove the support property. Let $c$ be a support constant of $\sigma$ and $E\in\mathcal{D}$ be $\tau$-semistable. Let $A_1,\dots,A_r$ be the HN factors of $E$ with respect to $\sigma$. By the definition of distance, we have
    \[\left|\log\dfrac{m_\sigma(E)}{m_\tau(E)}\right|\leqslant d(\sigma,\tau),\]
    and hence
    \[m_\sigma(E)\leqslant e^{d(\sigma,\tau)}m_\tau(E).\]
    Thus, we have
    \begin{align*}
        |Z_\tau(\cl E)|=m_\tau(E)&\geqslant e^{-d(\sigma,\tau)}m_\sigma(E)\\
        &=e^{-d(\sigma,\tau)}\sum_{j=1}^r|Z_\sigma(\cl A_j)|\\
        &\geqslant e^{-d(\sigma,\tau)}\sum_{j=1}^r c\|\cl A_j\|\\
        &\geqslant e^{-d(\sigma,\tau)}c\left\|\sum_{j=1}^r \cl A_j\right\|\\
        &=e^{-d(\sigma,\tau)}c\|\cl E\|.
    \end{align*}
    So $\tau$ has the support property with support constant $e^{-d(\sigma,\tau)}c$.

    Now let us show that $\tau$ satisfies the finite counting condition. What we need to show is that for any $\gamma\in\Gamma$ and $\psi\in\mathbb{R}$, the set
    \[\{E\in\Iso(\mathcal{P}_\tau(\psi))|\cl E=\gamma\}\]
    is finite. For any $0\neq E\in\mathcal{D}$, the definition of distance gives
    \[|\phi_\sigma^{\pm}(E)-\phi_\tau^{\pm}(E)|\leqslant d(\sigma,\tau),\left|\log\dfrac{m_\sigma(E)}{m_\tau(E)}\right|\leqslant d(\sigma,\tau).\]
    Hence,
    \begin{align*}
        m_\sigma(E)&\leqslant e^{d(\sigma,\tau)}m_\tau(E)\\
        \phi_\tau^-(E)-d(\sigma,\tau)\leqslant\phi_\sigma^-(E)&\leqslant\phi_\sigma^+(E)\leqslant\phi_\tau^+(E)+d(\sigma,\tau).
    \end{align*}
    Now let $E\in\mathcal{P}_\tau(\psi)$ with $\cl E=\gamma\neq0$. The two inequalities become
    \begin{align*}
        m_\sigma(E)&\leqslant e^{d(\sigma,\tau)}|Z_\tau(\gamma)|\\
        \psi-d(\sigma,\tau)\leqslant\phi_\sigma^-(E)&\leqslant\phi_\sigma^+(E)\leqslant\psi+d(\sigma,\tau).
    \end{align*}
    Apply Proposition \ref{mass and phase bounds finite} to
    \[R=e^{d(\sigma,\tau)}|Z_\tau(\gamma)|,I=[\psi-d(\sigma,\tau),\psi+d(\sigma,\tau)].\]
    Then the finiteness follows.
\end{proof}
\begin{proposition}\label{nearby}
    Let $\sigma=(Z_\sigma,\mathcal{P}_\sigma)\in\Stab_\cl(\mathcal{D})$ satisfy the support property and the finite counting condition with support constant $c$. For any $R>0$, there exists $\eta_R>0$, such that for any $\tau\in\Stab_\cl(\mathcal{D})$ with $d(\sigma,\tau)<\eta_R$, the following holds:

    \begin{enumerate}
        \item every $\tau$-semistable object $E$ with $m_\tau(E)\leqslant2R$ is $\sigma$-semistable;
        \item if $0\neq E\in\mathcal{P}_\sigma(\phi)$ satisfies $m_\sigma(E)\leqslant R$, then every HN factor of $E$ with respect to $\tau$ belongs to $\mathcal{P}_\sigma(\phi)$;
        \item $\|Z_\tau-Z_\sigma\|\leqslant\dfrac{c}4$, where $\|\cdot\|$ is the operator norm on the dual normed space $\Hom_\mathbb{Z}(\Gamma,\mathbb{C})=\Gamma_\mathbb{C}^*$.
    \end{enumerate}
\end{proposition}
\begin{proof}
    We give the proof by gradually picking smaller $\eta_R$ to adjust all the requirements.

    Firstly, take $0<\eta_0<\log2$. Then the condition $m_\tau(E)\leqslant2R$ implies that $m_\sigma(E)<4R$. By shifting by an even number, we can assume that $\phi_\tau(E)\in[0,2)$. So
    \[-\eta_0<\phi_\tau(E)-d(\sigma,\tau)\leqslant\phi_\sigma^-(E)\leqslant\phi_\sigma^+(E)\leqslant\phi_\tau(E)+d(\sigma,\tau)<2+\eta_0\]
    Applying Proposition \ref{mass and phase bounds finite} to $4R$ and $I=(-\eta,2+\eta),\eta>0$, we obtain that the set
    \[S_\eta=\{0\neq X\in\Iso(\mathcal{P}_\sigma(I))|m_\sigma(X)<4R\}\]
    is finite. We have the convention that the minimum of an empty set is $+\infty$. Now take
    \[0<\eta_1<\min\left(\dfrac12\min_{X\in S,\Delta_\sigma(X)>0}\Delta_\sigma(X),\eta_0\right),\]
    where $S=\{0\neq X\in\Iso(\mathcal{P}_\sigma(-1,3))|m_\sigma(X)<4R\}$, then every $\tau$-semistable object in $S_{\eta_1}$ is $\sigma$-semistable since its phase difference with respect to $\sigma$ cannot be positive. Hence, every $0<\eta\leqslant\eta_1$ satisfies the requirement (1).

    Now consider (2). Let $F_1,\dots,F_r$ be HN factors of $E$ with respect to $\tau$. Since $m_\sigma(E)\leqslant R$, we have $m_\tau(E)\leqslant 2R$, and hence $m_\tau(F_j)\leqslant m_\tau(E)<2R$, $j=1,\dots,r$. By (1), $F_j$ is $\sigma$-semistable, and
    \[\phi-\eta<\phi_\sigma(E)-d(\sigma,\tau)\leqslant\phi_\tau^-(E)\leqslant\phi_\tau(F_j)\leqslant\phi_\tau^+(E)\leqslant\phi_\sigma(E)+d(\sigma,\tau)<\phi+\eta,\]
    so $|\phi_\tau(F_j)-\phi|<\eta$. We also have $|\phi_\sigma(F_j)-\phi_\tau(F_j)|\leqslant d(\sigma,\tau)<\eta$, so
    \[|\phi_\sigma(F_j)-\phi|<2\eta.\]
    Now pick
    \[\eta_2=\min\left(\dfrac12\min_{X_1,X_2\in S\ \sigma-\text{semistable},\phi_\sigma(X_1)\neq\phi_\sigma(X_2)}|\phi_\sigma(X_1)-\phi_\sigma(X_2)|,\eta_1\right).\]
    Then when $0<\eta<\eta_2$, $F_j$ cannot be taken to have $\sigma$-phase different from $\phi$. Thus any $0<\eta<\eta_2$ satisfies (1) and (2).

    For (3), note that the map $\Stab_\cl(\mathcal{D})\rightarrow\Hom_{\mathbb{Z}}(\Gamma,\mathbb{C})$, sending each stability condition to its central charge, is continuous. Hence, we can pick $0<\eta_3<\eta_2$, such that (3) is satisfied. With all above, any $0<\eta_R<\eta_3$ satisfies all the requirements.
\end{proof}

\subsection{Stability structures of graded Lie algebras}
Let $\Lambda$ be a commutative cancellative monoid with zero that can be embedded into an abelian group. For example, $\mathbb{N}^k$ for some $k\in\mathbb{Z}_+$. Write $\Lambda_+=\Lambda-\{0\}$. Fix an additive \emph{height function} $\nu:\Lambda\rightarrow[0,+\infty)$, satisfying:
\begin{itemize}
    \item $\nu(\alpha)>0$ for any $\alpha\in\Lambda_+$;
    \item For any $R>0$, the sublevel set $\nu^{-1}([0,R])$ is finite.
\end{itemize}

Let
\[\mathfrak{g}=\prod_{\alpha\in\Lambda_+}\mathfrak{g}_\alpha\]
be a $\Lambda_+$-graded $\mathbb{Q}$-Lie algebra. For $R>0$, we define
\[F^{>R}\mathfrak{g}=\prod_{\nu(\alpha)>R}\mathfrak{g}_\alpha.\]
It is an ideal of $\mathfrak{g}$ by the additivity and positivity of $\nu$. Moreover the quotient $\mathfrak{g}/F^{>R}\mathfrak{g}$ is nilpotent since every sublevel set is finite. Thus we can take the pro-nilpotent group
\[G=\varprojlim_{R>0}\exp(\mathfrak{g}/F^{>R}\mathfrak{g}).\]

We denote by $H$ the strict upper half plane of $\mathbb{C}$, that is
\[H=\{z\in\mathbb{C}|\im z>0\}.\]
A rotation $e^{i\theta}H\subseteq\mathbb{C},\theta\in\mathbb{R}$ is called an \emph{open half plane} in the rest of this section.
\begin{definition}
    A \emph{central charge} on $\Lambda$ is an additive function $Z:\Lambda\rightarrow\mathbb{C}$ satisfying that $Z(\Lambda_+)\subseteq e^{i\theta}H$ for some $\theta\in\mathbb{R}$.
\end{definition}

Consider the multi-valued argument function $\Arg$. For any $\theta\in\mathbb{R}$, it has single-valued branches on the open half plane $e^{i\theta}H$. Fix one of those branches and denote it by $\arg$. It is worth warning that the symbol $\arg$ may not be the same as the usual branch of the argument function.

Let $Z:\Lambda\rightarrow\mathbb{C}$ be a central charge such that $Z(\Lambda_+)\subseteq e^{i\theta}H$ and $\arg$ be the branch of the argument function on $e^{i\theta}H$ as in the last paragraph. For $\alpha\in\Lambda_+$, we call
\[\phi_Z(\alpha)=\dfrac{1}{\pi}\arg Z(\alpha)\]
the \emph{phase} of $\alpha$ under $Z$. It is clear that the range of $\phi_Z$ lies in an open segment $(a,a+1)$ for some $a\in\mathbb{R}$. For $\phi\in(a,a+1)$, we set
\[\mathfrak{g}_\phi^Z=\prod_{\phi_Z(\alpha)=\phi}\mathfrak{g}_\alpha.\]
\begin{lemma}[Ordered factorization, {\cite[Section 2.2, Theorem 2 and Remark 7]{KS2008}}]\label{ordered decomposition}
    For any $g\in G$, there exists a unique family $L_\phi^Z\in\mathfrak{g}_\phi^Z,\phi\in (a,a+1)$, such that
    \[g=\prod_{\phi\in(a,a+1)}^\leftarrow\exp L_\phi^Z,\]
    where the product on the right-hand side is finite at any height cutoff.
\end{lemma}
An equivalent statement of Lemma \ref{ordered decomposition} is that the map
\begin{align*}
    \Phi_Z:\mathfrak{g}=\prod_{\alpha\in\Lambda_+}\mathfrak{g}_\alpha&\rightarrow G\\
    e=(e_{\alpha})&\mapsto\prod_{\phi\in(a,a+1)}^\leftarrow\exp\left(\sum_{\phi_Z(\alpha)=\phi}e_\alpha\right)
\end{align*}
is bijective. For two central charges $Z,W$ on $\Lambda$, we can define the \emph{transformation map}
\[T_{Z,W}=\Phi_W^{-1}\circ\Phi_Z:\mathfrak{g}\rightarrow G\rightarrow\mathfrak{g}.\]
It is easy to check that for three central charges $Z,W,V$ one has
\[T_{W,V}\circ T_{Z,W}=T_{Z,V},\]
and for any $Z$ one has
\[T_{Z,Z}=\id_{\mathfrak{g}}.\]

\begin{theorem}[{\cite[Section 2.3, Corollary 1]{KS2008}}]\label{universal Lie}
    Let $\alpha\in\Lambda_+$ and $Z,W$ be two central charges on $\Lambda$. There exists a Lie polynomial ${P_\alpha^{Z,W}}$ with coefficients in $\mathbb{Q}$ such that
    \[\pr_\alpha\circ T_{Z,W}(e)=P_\alpha^{Z,W}(e)\]
    for any $e\in\mathfrak{g}$.
\end{theorem}
\begin{proof}
    We first construct a free Lie algebra. For any $\beta\in\Lambda_+$, let $x_\beta$ be a formal symbol. Let $\mathfrak{l}$ be the completed free Lie algebra over $\mathbb{Q}$ generated by $x_\beta,\beta\in\Lambda_+$. It is a $\Lambda_+$-graded Lie algebra, and we can apply Lemma \ref{ordered decomposition} to it. Let $L$ be the corresponding pro-nilpotent group of $\mathfrak{l}$, and set
    \[l_Z=\prod_{\phi\in(a,a+1)}^\leftarrow\exp\left(\sum_{\phi_Z(\beta)=\phi}x_\beta\right).\]
    Take the ordered factorization of $l_Z$ with respect to $W$. Then its $\alpha$-entry must be an element in $\mathfrak{l}$, so it is a (possibly infinite) linear combination of brackets of $x_\beta$'s. 

    Now let us show that the linear combination must be finite. Note that a bracket of $x_{\beta_1},\dots,x_{\beta_s}$ contributes to an $\alpha$-entry of an element in $\mathfrak{l}$ only if $\beta_1+\dots+\beta_s=\alpha$. By the additivity of the height function, we have
    \[\nu(\alpha)=\nu(\beta_1)+\dots+\nu(\beta_s).\]
    By the positivity of $\nu$, we have $\nu(\beta_j)\leqslant\nu(\alpha),j=1,\dots,s$. So a generator $x_\beta$ occurs in the expression of $\pr_\alpha\circ\Phi_W^{-1}(l_Z)$ only if $\nu(\beta)\leqslant\nu(\alpha)$. Since the sublevel set $\nu^{-1}([0,\nu(\alpha)])$ is finite, only finitely many $x_\beta$ occur in the expression, and hence the linear combination must be finite.

    The above argument shows that there is a Lie polynomial $P^{Z,W}_\alpha$ such that
    \[\pr_\alpha\circ T_{Z,W}((x_\beta))=\pr_\alpha\circ\Phi_W^{-1}(l_Z)=P_\alpha^{Z,W}((x_\beta)).\]
    Let
    \begin{align*}
        \ev_e:\mathfrak{l}&\rightarrow\mathfrak{g}\\
        x_\beta&\mapsto e_\beta
    \end{align*}
    be the evaluation homomorphism at $e\in\mathfrak{g}$. One can easily show that $\ev_e$ commutes with any $T_{Z,W}$. Applying it to the last equation gives the result we need.
\end{proof}

\section{Obstruction to completion of derived Hall algebras}
Let $\mathcal{A}$ be an abelian category of finite length. Recall in \cite{BT2012} Bridgeland and Toledano Laredo provide a formula
\[\1_{\mathcal{A}}=\prod_{\phi\in(0,1]}^\leftarrow\Ss^Z_{\phi}\]
in the completion $\hat{\mathcal{H}}(\mathcal{A})$ of the Ringel-Hall algebra $\mathcal{H}(\mathcal{A})$. For $\mathcal{D}$ and its derived Hall algebra, we want to do the same. That is the goal of this section. However, it turns out that one cannot define such a completion of $\mathcal{H}(\mathcal{D})$. More precisely, for $U\subseteq\Iso(\mathcal{D})$, we write $\mathscr{F}(U)$ for the set of $\mathbb{Q}$-valued functions on $\Iso(\mathcal{D})$ supported on $U$. Then the multiplication (\ref{eq:Hall}) cannot be extended continuously to the whole $\mathscr{F}(\Iso(\mathcal{D}))$. We will explain the details in this section.

\subsection{Topologies on function spaces}
Let $U\subseteq\Iso(\mathcal{D})$. We will define a topology on the function space $\mathscr{F}(U)$. Note that
\[\mathscr{F}(U)=\prod_{X\in U}\mathbb{Q}u_X.\]
On each $\mathbb{Q}u_X$, we pick the discrete topology, and we take the product topology on $\mathscr{F}(U)$.

It is worth warning that the product of discrete topologies need not be discrete. When $U$ is infinite, the topology we take is quite different from the discrete topology. In fact, by the definition of product topology, any open set in $\mathscr{F}(U)$ can be written as a union of subsets of the form
\[\prod_{X\in U-\{X_1,\dots,X_k\}}\mathbb{Q}u_X\times A_1u_{X_1}\times\dots\times A_ku_{X_k},\]
where $X_1,\dots,X_k\in U$ and $A_1,\dots,A_k\subseteq\mathbb{Q}$. In particular, any one-point subset of $\mathscr{F}(U)$ is not open when $U$ is infinite. In fact, by the properties of product topologies, we have the following lemma:
\begin{lemma}\label{neighbourhood basis}
    For any $f\in\mathscr{F}(U)$, it has a neighbourhood basis $\{N_{f,F}^U\}$, where $F$ runs over all finite subsets of $U$ and
\[N_{f,F}^U=\{g\in\mathscr{F}(U)|g|_F=f|_F\}.\]
\end{lemma}
When $f=0\in\mathscr{F}(U)$, we write $N_F^U=N_{0,F}^U$ for short. Then it is easy to check that $N_{f,F}^U=f+N_F^U$ for any $f\in\mathscr{F}(U)$.

We denote by $\mathscr{F}_{\fin}(U)$ the subset of $\mathscr{F}(U)$ consisting of functions that have finite support. In fact,
\[\mathscr{F}_{\fin}(U)=\bigoplus_{X\in U}\mathbb{Q}u_X.\]
Then we have
\begin{lemma}\label{density}
    $\mathscr{F}_{\fin}(U)$ is dense in $\mathscr{F}(U)$.
\end{lemma}
\begin{proof}
    It suffices to show that, for any $f\in\mathscr{F}(U)$ and any neighbourhood $N$ of $f$, the intersection $N\cap\mathscr{F}_{\fin}(U)$ is not empty. By Lemma \ref{neighbourhood basis}, one only needs to check for $N=N^U_{f,F}$ for any finite subset $F$ of $U$. Consider the function
    \[f_F=\sum_{X\in F}f(X)u_X.\]
    This function is finitely supported and agrees with $f$ on $F$, and hence belongs to $N^U_{f,F}\cap\mathscr{F}_{\fin}(U)$.
\end{proof}

We also have a description of convergence in $\mathscr{F}(U)$:
\begin{lemma}\label{convergence}
    Let $\{f_n\}$ be a sequence in $\mathscr{F}(U)$. $\{f_n\}$ converges to $f\in\mathscr{F}(U)$ if and only if for any finite subset $F\subseteq U$, there exists an $N>0$ such that $f|_F=f_n|_F$ for any $n>N$.
\end{lemma}
It is a direct consequence of Lemma \ref{neighbourhood basis} and we omit the proof.
\subsection{Obstruction to extending multiplication}
In this section, we give a proof that we cannot extend the multiplication (\ref{eq:Hall}) continuously to $\mathscr{F}(\Iso(\mathcal{D}))$, in the sense of the topology defined in the last subsection.
\begin{definition}
    Let $U,V\subseteq\Iso(\mathcal{D})$ and $X\in\Iso(\mathcal{D})$. We set
    \[R_X(U,V)=\{(A,B)\in U\times V|F_{AB}^X\neq0\}.\]
    We say the pair $(U,V)$ is \emph{admissible} if $R_X(U,V)$ is finite for any $X\in\Iso(\mathcal{D})$.
\end{definition}
The following theorem is the criterion of extension.
\begin{theorem}[Continuous extension criterion]\label{continuous extension criterion}
    The multiplication
    \[*:\mathscr{F}_{\fin}(U)\times\mathscr{F}_{\fin}(V)\rightarrow\mathscr{F}_{\fin}(\Iso(\mathcal{D}))=\mathcal{H}(\mathcal{D})\]
    extends to a jointly continuous bilinear map
    \[\mathscr{F}(U)\times\mathscr{F}(V)\rightarrow\mathscr{F}(\Iso(\mathcal{D}))\]
    if and only if $(U,V)$ is admissible. In this case, the extension is unique and is given by
    \[(f*g)(X)=\sum_{(A,B)\in R_X(U,V)}F_{AB}^Xf(A)g(B).\]
\end{theorem}
\begin{proof}
    Let us first show the ``only if'' part. Let
    \[m:\mathscr{F}(U)\times\mathscr{F}(V)\rightarrow\mathscr{F}(\Iso(\mathcal{D}))\]
    be an extension of $*$. Given $X\in\Iso(\mathcal{D})$, we define a function
    \begin{align*}
        B_X:\mathscr{F}(U)\times\mathscr{F}(V)&\rightarrow\mathbb{Q}\\
        (f,g)&\mapsto m(f,g)(X).
    \end{align*}
    Since $m$ is jointly continuous, $B_X=\pr_X\circ m$ is also jointly continuous, where the topology of the codomain $\mathbb{Q}$ is discrete. Pick the open subset $\{0\}$ in $\mathbb{Q}$. Then $B_X^{-1}(0)$ is open in $\mathscr{F}(U)\times\mathscr{F}(V)$. By Lemma \ref{neighbourhood basis}, we can choose finite subsets $F\subseteq U,G\subseteq V$ such that
    \[N_F^U\times N_G^V\subseteq B_X^{-1}(0).\]

    If $A\in U-F,B\in V-G$, then $u_A\in N_F^U$ and $u_B\in N_G^V$. By the hypothesis of $F$ and $G$, we have
    \[0=B_X(u_A,u_B)=m(u_A,u_B)(X)=(u_A*u_B)(X)=F_{AB}^X.\]
    As a consequence,
    \[R_X(U,V)\subseteq(F\times V)\cup (U\times G).\]
    Let $R_F=R_X(U,V)\cap(F\times V), R_G=R_X(U,V)\cap(U\times G)$. Then a simple calculation of sets shows that $R_X(U,V)=R_F\cup R_G$. We will show that $R_F$ is finite next. Then the same argument gives that $R_G$ is finite, and hence $R_X(U,V)$ is finite.

    For any $A\in F$, we set
    \[R_A=\{(Y,B)\in R_F|Y=A\},\]
    so we have a partition of $R_F$:
    \[R_F=\bigsqcup_{A\in F}R_A.\]
    So it suffices to show that each $R_A$ is finite. Consider the map
    \begin{align*}
        \ell_A:\mathscr{F}(V)&\rightarrow\mathbb{Q}\\
        g&\mapsto B_X(u_A,g).
    \end{align*}
    It is continuous and linear, so its kernel is open. By Lemma \ref{neighbourhood basis}, there is a finite subset $G_A\subseteq V$ such that $N_{G_A}^V\subseteq\Ker\ell_A$. For $B\in V-G_A$, we have $u_B\in N_{G_A}^V$, so
    \[0=\ell_A(u_B)=B_X(u_A,u_B)=F_{AB}^X.\]
    Thus, we know that $R_A\subseteq \{A\}\times G_A$ is finite. This proves the ``only if'' part.

    For the ``if'' part, note that
    \[m(f,g)(X)=\sum_{(A,B)\in R_X(U,V)}F_{AB}^Xf(A)g(B)\]
    is well-defined. One can easily check that $m$ satisfies the joint continuity condition and agrees with $*$ on $\mathscr{F}_{\fin}(U)\times\mathscr{F}_{\fin}(V)$.

    Finally, by Lemma \ref{density}, $\mathscr{F}_{\fin}(U)\times\mathscr{F}_{\fin}(V)$ is dense in $\mathscr{F}(U)\times\mathscr{F}(V)$, and the uniqueness follows.
\end{proof}
Now we can state and prove the obstruction we mentioned in the beginning of this section.
\begin{theorem}\label{obstruction}
The multiplication (\ref{eq:Hall}) cannot be extended continuously to the whole space $\mathscr{F}(\Iso(\mathcal{D}))$.
\end{theorem}
\begin{proof}
    Let $0\neq S\in\mathcal{D}$. When $n\neq m$, we have $S[n]\not\cong S[m]$. In fact, if $S[n]\cong S[m]$, then $S\cong S[d]$ for some $d>0$ and hence
    \[\sum_{i\geq0}\dim_k\Hom_{\mathcal{D}}(S[i],S)\geqslant\sum_{j\geq0}\dim_k\Hom_{\mathcal{D}}(S[jd],S)=\infty,\]
    which contradicts the left homological finiteness condition (\ref{eq:leftfinite}). As a consequence, we deduce that the sequence $\{u_{S[n]}\}$ converges to $0$ in $\mathscr{F}(\Iso(\mathcal{D}))$, by using Lemma \ref{convergence}.
    
     If there is a jointly continuous extension of $*$ on the whole space $\mathscr{F}(\Iso(\mathcal{D}))$, then we must have
    \[\lim_{n\rightarrow\infty}F_{S[n]S[n+1]}^0=\lim_{n\rightarrow\infty}B_0(u_{S[n]},u_{S[n+1]})=0.\]
    However, let $A=S[n]$ and $B=S[n+1]=A[1]$. Then
    \[F_{AB}^0=\dfrac{|\Hom_{\mathcal{D}}(0,B)_{A[1]}|}{a_B}\dfrac{\{0,B\}}{\{B,B\}}=\dfrac{1}{a_B\{B,B\}}=\dfrac{1}{a_S\{S,S\}},\]
    which is positive and independent of $n$. Hence
    \[\lim_{n\rightarrow\infty}F_{S[n]S[n+1]}^0=\dfrac{1}{a_S\{S,S\}}>0,\]
    a contradiction. That shows that the global extension does not exist.
\end{proof}
\subsection{The degreewise completion}
In the original version \cite{GX2025} of this work, we considered the degreewise completion
\[\hat{\mathcal{H}}(\mathcal{D})=\prod_{\gamma\in K_0(\mathcal{D})}\mathcal{H}(\mathcal{D})_\gamma.\]
However, that is not the correct choice by the following proposition.
\begin{proposition}
    Let $\1_{\mathcal{D}}\in\mathscr{F}(\Iso(\mathcal{D}))$ be the function taking the value $1$ on every isomorphism class in $\Iso(\mathcal{D})$. Then
    \[\1_{\mathcal{D}}\notin\hat{\mathcal{H}}(\mathcal{D}).\]
\end{proposition}
\begin{proof}
Let $0\neq S\in\mathcal{D}$. Then for any $n\in\mathbb{Z}$, we have
\[\udv S[2n]=\udv S.\]
Since $S\neq 0$, the beginning part of the proof of Theorem \ref{obstruction} shows that $S[2n]\not\cong S[2m]$ for $n\neq m$. So the projection of $\1_{\mathcal{D}}$ on the subspace $\mathcal{H}(\mathcal{D})_{\udv S}$ is not finitely supported. That is a contradiction.
\end{proof}
\subsection{The descending product}
Although we cannot work on the degreewise completion of $\mathcal{H}(\mathcal{D})$, we can still state \cite[Theorem 4.19]{GX2025} in $\mathscr{F}(\Iso(\mathcal{D}))$, but we have to understand it in a correct way.

Let us do some preparation first. Let $\sigma=(Z,\mathcal{P})$ be a stability condition on $\mathcal{D}$. We have the following lemma.
\begin{lemma}\label{segment admissible}
    Let $I_1>I_2$ be segments in $\mathbb{R}$. Then the pair $(\Iso(\mathcal{P}(I_1)),\Iso(\mathcal{P}(I_2)))$ is admissible.
\end{lemma}
\begin{proof}
    For any $X\in\Iso(\mathcal{D})$, we are going to show that the set
    \[R_X(\Iso(\mathcal{P}(I_1)),\Iso(\mathcal{P}(I_2)))
    \] is finite. We will write $U_1=\Iso(\mathcal{P}(I_1)),U_2=\Iso(\mathcal{P}(I_2))$ for convenience.

    If $A\in U_1,B\in U_2$ satisfy the condition $F_{AB}^X\neq0$, then there exists a triangle
    \[A\rightarrow X\rightarrow B\rightarrow A[1]\]
    in $\mathcal{D}$. Since $I_1>I_2$, there exists a $c\in\mathbb{R}$ such that $I_1\subseteq(c,+\infty)$ and $I_2\subseteq(-\infty,c]$ or $I_1\subseteq[c,+\infty)$ and $I_2\subseteq(-\infty,c)$. Without loss of generality, we may assume the first condition holds. Let
    \[X_{>c}\rightarrow X\rightarrow X_{\leqslant c}\rightarrow X_{>c}[1]\]
    be the phase truncation of $X$. By the uniqueness part of Proposition \ref{phase truncation}, we have $A\cong X_{>c}$ and $B\cong X_{\leqslant c}$. Hence, $R_X(U_1,U_2)$ is a one-point set if it is not empty, which is obviously finite.
\end{proof}
In fact, the derived Hall number in this case can be calculated.
\begin{lemma}\label{truncation hall number 1}
    Let $X_{>c}\rightarrow X\rightarrow X_{\leqslant c}\rightarrow X_{>c}[1]$ be the phase truncation of $X$ at $c\in\mathbb{R}$. Then the derived Hall number
    \[F_{X_{>c}X_{\leqslant c}}^X=1.\]
\end{lemma}
\begin{proof}
    By the definition of derived Hall number we have
    \[F_{X_{>c}X_{\leqslant c}}^X=\dfrac{|\Hom_{\mathcal{D}}(X,X_{\leqslant c})_{X_{>c}[1]}|}{a_{X_{\leqslant c}}}\dfrac{\{X,X_{\leqslant c}\}}{\{X_{\leqslant c},X_{\leqslant c}\}}.\]
    By Proposition \ref{phase truncation}, the phase truncation is unique up to unique isomorphism, so
    \[|\Hom_{\mathcal{D}}(X,X_{\leqslant c})_{X_{>c}[1]}|=a_{X_{\leqslant c}}.\]
    Now apply the cohomological functor $\Hom_{\mathcal{D}}(-,X_{\leqslant c})$ to the phase truncation of $X$ at $c$. We obtain the long exact sequence
    \begin{align*}
        \dots\rightarrow\Hom(X_{>c}[n+1],X_{\leqslant c})&\rightarrow\Hom(X_{\leqslant c}[n],X_{\leqslant c})\xrightarrow{f_n}\Hom(X[n],X_{\leqslant c})\\
        &\rightarrow \Hom(X_{>c}[n],X_{\leqslant c})\rightarrow\cdots.
    \end{align*}
    When $n>0$, $\Hom(X_{>c}[n],X_{\leqslant c})=0$ by the definition of a Bridgeland stability condition, and hence $f_n$ is an isomorphism. That tells us that $\{X,X_{\leqslant c}\}=\{X_{\leqslant c},X_{\leqslant c}\}$.
\end{proof}
\begin{corollary}
    Let $I_1>I_2>\dots>I_r$ be segments in $\mathbb{R}$ and $U_j=\Iso(\mathcal{P}(I_j))$ and $f_j\in\mathscr{F}(U_j)$ for $j=1,\dots,r$. Then the product $f_1*f_2*\dots*f_r$ is well defined in $\mathscr{F}(\Iso(\mathcal{P}(I_1)*\dots*\mathcal{P}(I_r)))$ with the explicit expression:
    \[(f_1*f_2*\dots*f_r)(X)=\prod_{j=1}^r f_j(X_j),\]
    where $X_j\in\mathcal{P}(I_j)$ is the canonical block of $X$ at $I_j$ for each $j=1,\dots,r$. 
\end{corollary}
\begin{proof}
    Let us do induction on $r$. When $r=1$ there is nothing to prove. When $r=2$ it is a natural consequence of Lemma \ref{segment admissible} and Lemma \ref{truncation hall number 1}. Now suppose that the corollary holds for $1,2,\dots, r-1$. Let $I\subseteq\mathbb{R}$ be a segment such that $I_1,I_2,\dots,I_{r-1}\subseteq I$ and $I\cap I_r=\varnothing$. Set $U=\Iso(\mathcal{P}(I))$. By the induction hypothesis, the multiplication (\ref{eq:Hall}) continuously extends to
    \begin{align*}
        \mathscr{F}(U_1)\times\mathscr{F}(U_2)\times\dots\times\mathscr{F}(U_{r-1})&\rightarrow\mathscr{F}(\Iso(\mathcal{P}(I_1)*\dots*\mathcal{P}(I_{r-1})))\\
        (f_1,f_2,\dots,f_{r-1})&\mapsto\left(X\mapsto\prod_{j=1}^{r-1}f_j(X_j)\right).
    \end{align*}
    The image of the extension is in fact included in $\mathscr{F}(U)$. Now applying the $r=2$ case to $I$ and $I_r$, one obtains that the extension exists and the explicit expression is
    \begin{align*}
        (f_1*f_2*\dots*f_r)(X)&=((f_1*f_2*\dots*f_{r-1})*f_r)(X)\\
        &=(f_1*f_2*\dots*f_{r-1})(X_I)f_r(X_r)\\
        &=\prod_{j=1}^{r-1}f_j(X_j)\cdot f_r(X_r)\\
        &=\prod_{j=1}^{r}f_j(X_j).
    \end{align*}
\end{proof}
Now we can state the descending product formula.
\begin{theorem}[Descending product]\label{descending product}
    The formula
    \begin{equation}\label{eq:descending product}
        \1_{\mathcal{D}}=\prod_{\phi\in\mathbb{R}}^\leftarrow\1_{\mathcal{P}(\phi)}
    \end{equation}
    holds coefficientwise. That is, for any $X\in\Iso(\mathcal{D})$, the coefficients of $u_X$ in both sides of (\ref{eq:descending product}) are well defined and equal.
\end{theorem}
\begin{proof}
    If $X=0$, then the right-hand side is a product of $1$'s, and hence equal to $1$. If $X\neq 0$, let $X_1,\dots, X_r$ be the HN factors of $X$ with $\phi(X_1)>\dots>\phi(X_r)$. So the right-hand side is
    \[\left(\prod_{\phi\in\mathbb{R}}^\leftarrow\1_{\mathcal{P}(\phi)}\right)(X)=(\1_{\mathcal{P}(\phi(X_1))}*\dots*\1_{\mathcal{P}(\phi(X_r))})(X).\]
    Take segments $I_1>\dots>I_r$ such that $\phi(X_j)\in I_j,j=1,\dots,r$. Then by Lemma \ref{truncation hall number 1},
    \begin{align*}
        \left(\prod_{\phi\in\mathbb{R}}^\leftarrow\1_{\mathcal{P}(\phi)}\right)(X)&=(\1_{\mathcal{P}(\phi(X_1))}*\dots*\1_{\mathcal{P}(\phi(X_r))})(X)\\
        &=\1_{\mathcal{P}(\phi(X_1))}(X_1)\cdots\1_{\mathcal{P}(\phi(X_r))}(X_r)\\
        &=1\\
        &=\1_{\mathcal{D}}(X),
    \end{align*}
    where we use the fact that the canonical HN block of $X$ in $I_j$ is $X_j$, $j=1,\dots,r$.
\end{proof}
\begin{remark}
    The right-hand side of (\ref{eq:descending product}) can be defined as a limit in $\mathscr{F}(\Iso(\mathcal{D}))$. To make this precise, for $F=\{\phi_1,\dots,\phi_r\}\subseteq\mathbb{R}$ where $\phi_1>\dots>\phi_r$, set
    \[p_F=\1_{\mathcal{P}(\phi_1)}*\cdots*\1_{\mathcal{P}(\phi_r)},\]
    with $p_{\varnothing}=u_0$. These finite ordered products are well defined by the preceding results. For $X\in\Iso(\mathcal{D})$, let $\Phi_X$ be the finite set of phases of its nonzero HN factors, with $\Phi_0=\varnothing$. By the uniqueness of HN filtrations and Lemma \ref{truncation hall number 1}, we have
    \[p_F(X)=\begin{cases}
            1, & \Phi_X\subseteq F,\\
            0, & \Phi_X\not\subseteq F.
        \end{cases}\]
    It follows that the net $(p_F)_F$, where $F$ runs over the finite subsets of $\mathbb{R}$ and is ordered by inclusions, converges to $\1_{\mathcal{D}}$ in the product topology:
    \[\lim_F p_F=\1_{\mathcal{D}}\in\mathscr{F}(\Iso(\mathcal{D})).\]
    With this interpretation of the infinite ordered product, (\ref{eq:descending product}) is an equality of functions.

    This construction uses only the finite ordered products of phase characteristic functions and their limit. It does not define a jointly continuous extension of Hall multiplication for arbitrary pairs of functions in $\mathscr{F}(\Iso(\mathcal{D}))$. Hence it is compatible with Theorem \ref{obstruction}.
\end{remark}

\section{Strict sectors and finite numerical Joyce transport}
In \cite{BT2012}, Bridgeland and Toledano Laredo provided the diagram (\ref{BT diagram}) to characterize the relationships between three types of elements. In the original version \cite{GX2025} of this work, we tried to recover it in $\mathcal{D}$. However, in fact, we cannot recover it in the whole $\mathcal{D}$, due to Theorem \ref{obstruction}. In this section, we will recover the diagram on segments.
\subsection{Strict sectors and characteristic functions}
For a \emph{strict sector}, we mean a segment $I\subseteq[b,d]\subseteq\mathbb{R}$ with $d-b<1$. We allow $I$ to be a single point $\{\phi\}$. Let $\sigma\in\Stab_\cl(\mathcal{D})$ satisfy the support property and the finite counting condition. We set
\[\mathcal{C}=\mathcal{P}_\sigma(I),S_I=\{\cl E|E\in\mathcal{C}\}.\]
Let $\mathcal{H}(\mathcal{C})$ be the subalgebra of $\mathcal{H}(\mathcal{D})$ generated by $u_X,X\in\mathcal{C}$. By the extension closedness of  $\mathcal{C}$, it is spanned by $u_X,X\in\mathcal{C}$ as a $\mathbb{Q}$-vector space. It has a natural grading
\[\mathcal{H}(\mathcal{C})=\bigoplus_{\gamma\in S_I}\mathcal{H}(\mathcal{C})_\gamma.\]
\begin{lemma}\label{coincide Ringel-Hall}
    $\mathcal{H}(\mathcal{C})$ coincides with the Ringel-Hall algebra of the exact category $\mathcal{C}$.
\end{lemma}
\begin{proof}
    Let $a\in\mathbb{R}$ such that $I\subseteq(a,a+1]$. Then the category $\mathcal{C}$ is included in the rotated heart $\mathcal{A}=\mathcal{P}_\sigma(a,a+1]\subseteq\mathcal{D}$, which is an abelian category. Let $\mathcal{H}(\mathcal{A})\subseteq\mathcal{H}(\mathcal{D})$ be the subalgebra generated by $u_X,X\in\mathcal{A}$. By a result by To\"en \cite[Proposition 6.1]{Toe2006}, $\mathcal{H}(\mathcal{A})$ coincides with the Ringel-Hall algebra of $\mathcal{A}$.
    
    Now $\mathcal{H}(\mathcal{C})$ is a subalgebra of $\mathcal{H}(\mathcal{A})$. Note that the exact structure on $\mathcal{C}$ is inherited from $\mathcal{A}$. So $\mathcal{H}(\mathcal{C})$ coincides with the Ringel-Hall algebra of $\mathcal{C}$.
\end{proof}
\begin{lemma}\label{sector degree finite}
    Each $\mathcal{H}(\mathcal{C})_\gamma$ is finite dimensional. Equivalently, for any $\gamma\in S_I$, the set
    \[\mathcal{C}_\gamma=\{E\in\Iso(\mathcal{C})|\cl E=\gamma\}\]
    is finite.
\end{lemma}
\begin{proof}
    We define
    \begin{align*}
        h:S_I&\rightarrow\mathbb{R}\\
        \gamma&\mapsto\re(e^{-i\pi\frac{b+d}{2}}Z_\sigma(\gamma)).
    \end{align*}
    The strictness of $I$ ensures that the range of $h$ lies in $[0,+\infty)$ and $h(\gamma)=0$ if and only if $\gamma=0$. Let $0\neq E\in\mathcal{C}$ with HN factors $A_1,\dots,A_r$. The support property gives
    \[h(\cl E)\geqslant a\sum_{j=1}^r|Z_\sigma(\cl A_j)|\geqslant ac\sum_{j=1}^r\|\cl A_j\|,\]
    where
    \[a=\cos\dfrac{\pi(d-b)}{2}>0.\]
    This inequality, together with the Hom-finiteness of $\mathcal{D}$ and the finite counting condition of $\sigma$, implies that the sublevel set $h^{-1}([0,R))$ is finite for any $R>0$. Thus, as a subset of $\{E\in\Iso(\mathcal{C})|\cl E\in h^{-1}([0,h(\gamma)+1))\}$, the set $\mathcal{C}_\gamma$ must be finite.
\end{proof}
As a consequence, we can define the characteristic function $\kappa^\sigma_{I,\gamma}\in\mathcal{H}(\mathcal{C})$ supported on the finite set $\mathcal{C}_\gamma$ for $\gamma\in S_I$. Set
\[\mathcal{P}_\sigma(\phi)_\gamma=\Iso(\mathcal{P}_\sigma(\phi))\cap\mathcal{C}_\gamma,\gamma\in S_I,\phi\in I.\]
It is also a finite set, and we denote the characteristic function supported on it by $\delta^\sigma_{\phi,\gamma}\in\mathcal{H}(\mathcal{C})$ when $\gamma\neq0$. It is clear that $\delta^\sigma_{\phi,\gamma}\neq0$ only if $\phi=\phi_\sigma(\gamma)$.
\begin{theorem}[HN identity and Reineke inversion]
    For any $0\neq\gamma\in S_I$, we have the \emph{HN identity}:
    \begin{equation}\label{HN identity in sector}
        \kappa^\sigma_{I,\gamma}=\sum_{\begin{subarray}{c}
         \\  {\gamma_1+\gamma_2+\dots+\gamma_r=\gamma},\gamma_i\neq0\\{\phi_\sigma(\gamma_1)>\phi_\sigma(\gamma_2)>\dots>\phi_\sigma(\gamma_r)\in I}
    \end{subarray}}\delta^\sigma_{\phi_\sigma(\gamma_1),\gamma_1}*\delta^\sigma_{\phi_\sigma(\gamma_2),\gamma_2}*\dots*\delta^\sigma_{\phi_\sigma(\gamma_r),\gamma_r},
    \end{equation}
    and its inverted version, called the \emph{Reineke inversion}:
    \begin{equation}\label{Reineke inversion in sector}
        \delta^\sigma_{\phi_\sigma(\gamma),\gamma}=\sum_{\begin{subarray}{c}
        \gamma_1+\gamma_2+\dots+\gamma_r=\gamma,\gamma_i\neq0\\\phi_\sigma(\gamma_1+\dots+\gamma_i)>\phi_\sigma(\gamma),1\leqslant i<r\\
    \end{subarray}}(-1)^{r-1}\,\kappa^\sigma_{I,\gamma_1}*\dots*\kappa^\sigma_{I,\gamma_r}.
    \end{equation}
\end{theorem}
\begin{proof}
    Since $\mathcal{C}$ is extension closed, every $0\neq E\in\mathcal{C}$ has a unique HN filtration with each HN factor in $\mathcal{C}$. The corresponding derived Hall numbers are all $1$ by Lemma \ref{truncation hall number 1}, and (\ref{HN identity in sector}) follows.

    Let $a\in\mathbb{R}$ such that $I\subseteq(a,a+1]$. Then the category $\mathcal{C}$ is included in the rotated heart $\mathcal{A}=\mathcal{P}_\sigma(a,a+1]\subseteq\mathcal{D}$. So the proof of (\ref{Reineke inversion in sector}) from (\ref{HN identity in sector}) is the same as in an abelian category, by Lemma \ref{coincide Ringel-Hall}, and we refer the reader to \cite[Theorem 5.4]{BT2012} and \cite[Theorem 5.1]{Rei2003}.
\end{proof}
Let $\phi\in I$ and $\gamma\in S_I$ with $\phi_\sigma(\gamma)=\phi$. We define \emph{Joyce's element} $\epsilon^\sigma_{\phi,\gamma}$ to be the finite sum
\begin{equation}\label{Joyce's element}
    \epsilon^\sigma_{\phi,\gamma}=\sum_{\begin{subarray}{c}
    \\ \gamma_1+\dots+\gamma_r=\gamma\\
   \phi_\sigma(\gamma_1)=\dots=\phi_\sigma(\gamma_r)=\phi_\sigma(\gamma) 
\end{subarray}}\frac{(-1)^{r-1}}{r}\delta^\sigma_{\phi,\gamma_1}*\dots*\delta^\sigma_{\phi,\gamma_r}.
\end{equation}
This expression is, in fact, a finite cutoff of a logarithmic series, so it has an inversion in the form of an exponential series:
\begin{equation}
    \label{Joyce inversion}
    \delta^\sigma_{\phi,\gamma}=\sum_{\begin{subarray}{c}
    \\ \gamma_1+\dots+\gamma_r=\gamma\\
   \phi_\sigma(\gamma_1)=\dots=\phi_\sigma(\gamma_r)=\phi_\sigma(\gamma) 
\end{subarray}}\frac{1}{r!}\epsilon^\sigma_{\phi,\gamma_1}*\dots*\epsilon^\sigma_{\phi,\gamma_r}.
\end{equation}
It is easy to see that for any $\phi$ and $\gamma$, the functions $\epsilon^\sigma_{\phi,\gamma}$ and $\delta^\sigma_{\phi,\gamma}$ in $\mathcal{H}(\mathcal{D})$ are independent of the choice of the sector $I\ni\phi$.

So far, we have reconstructed the diagram (\ref{BT diagram}) by Bridgeland and Toledano Laredo in the strict sector:
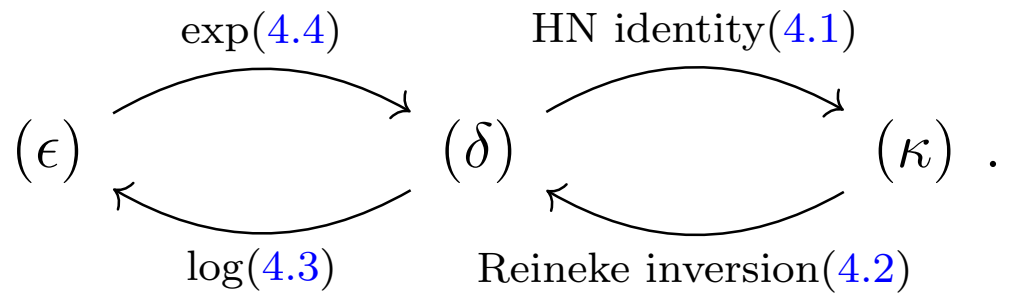
\begin{center}
    \begin{tikzcd}
(\epsilon) \arrow[rr, "\exp(\ref{Joyce inversion})", bend left] &  & (\delta) \arrow[ll, "\log(\ref{Joyce's element})", bend left] \arrow[rr, "\text{HN identity}(\ref{HN identity in sector})", bend left] &  & (\kappa) \arrow[ll, "\text{Reineke inversion}(\ref{Reineke inversion in sector})", bend left]
\end{tikzcd}.
\end{center}
\subsection{Degreewise completion in a strict sector}
In this subsection, we are going to construct a degreewise completion
\[\hat{\mathcal{H}}(\mathcal{C})=\prod_{\gamma\in S_I}\mathcal{H}(\mathcal{C})_\gamma.\]
Since each $\mathcal{H}(\mathcal{C})_\gamma$ is finite-dimensional, by Lemma \ref{sector degree finite}, $\hat{\mathcal{H}}(\mathcal{C})$ is the same space as $\mathscr{F}(\Iso(\mathcal{C}))$.
\begin{proposition}\label{sector completion}
    The multiplication (\ref{eq:Hall}) extends to $\hat{\mathcal{H}}(\mathcal{C})$ continuously via the finite degreewise expression
    \begin{equation}\label{eq:degreewise}
        (f*g)_\gamma=\sum_{\alpha+\beta=\gamma}f_\alpha*g_\beta.
    \end{equation}
\end{proposition}
\begin{proof}
    By Theorem \ref{continuous extension criterion}, we only need to check that for any $X\in\mathcal{C}$, the set
    \[R_X(\Iso(\mathcal{C}),\Iso(\mathcal{C}))\]
    is finite.
    
    For any $(A,B)\in R_X(\Iso(\mathcal{C}),\Iso(\mathcal{C}))$, since $F_{AB}^X\neq0$, we have a triangle
    \[A\rightarrow X\rightarrow B\rightarrow A[1].\]
    So
    \[\cl X=\cl A+\cl B.\]
    Recall the height function $h$ on $S_I$ defined in the proof of Lemma \ref{sector degree finite}. We obtain
    \[h(\cl X)= h(\cl A)+h(\cl B),\]
    and hence $h(\cl A),h(\cl B)\leqslant h(\cl X)$. Since the sublevel set $h^{-1}([0,h(\cl X)+1))$ is finite, together with Lemma \ref{sector degree finite}, both $A$ and $B$ have only finitely many choices. So $R_X(\Iso(\mathcal{C}),\Iso(\mathcal{C}))$ is finite.

    For the expression (\ref{eq:degreewise}), we only need to show that the right-hand side is a finite sum. This follows from the fact that $\alpha+\beta=\gamma$ implies $h(\alpha)+h(\beta)=h(\gamma)$ and the finiteness of sublevel sets of $h$.
\end{proof}
\subsection{Completed Hall algebras of phase categories}
In this subsection, we study the phase categories. Recall that by a \emph{phase category}, we mean a subcategory $\mathcal{P}_\sigma(\phi)\subseteq\mathcal{D}$ for some $\phi\in\mathbb{R}$. It is abelian and of finite length under our assumption.

Fix a phase category $\mathcal{B}=\mathcal{P}_\sigma(\phi),\phi\in\mathbb{R}$. Set
\[\Lambda_\phi=S_{\{\phi\}}=\{\cl E|E\in\mathcal{B}\}\subseteq\Gamma,\Lambda_\phi^+=\Lambda_\phi-\{0\}.\]
It is easy to check that $\Lambda_\phi$ is a commutative cancellative monoid with zero.
\begin{lemma}
    The function
    \begin{align*}
        h_\phi:\Lambda_\phi&\rightarrow[0,+\infty)\\
        \beta&\mapsto|Z_\sigma(\beta)|
    \end{align*}
    is a height function on $\Lambda_\phi$. Moreover, for any $0\neq E\in\mathcal{B}$, we have
    \[m_\sigma(E)=h_\phi(\cl E).\]
\end{lemma}
\begin{proof}
    Take
    \[I=\{\phi\},b=\phi-\dfrac14,d=\phi+\dfrac14,\]
    so $I\subseteq[b,d]$ and $\Lambda_\phi=S_I$. Now the height function defined in the proof of Lemma \ref{sector degree finite} is
    \[h(\gamma)=\re(e^{-i\pi\frac{b+d}{2}}Z_\sigma(\gamma))=\re(e^{-i\pi\phi}Z_\sigma(\gamma))=|Z_{\sigma}(\gamma)|,\]
    where the last identity follows from the fact that $Z_{\sigma}(\gamma)$ has phase $\phi$. So $h_\phi=h$ is a height function.

    Now let $0\neq E\in\mathcal{B}$. Since $E$ is $\sigma$-semistable, by the definition of mass, 
    \[m_\sigma(E)=|Z_\sigma(\cl E)|=h_\phi(\cl E).\]
\end{proof}

Let $\mathcal{H}(\mathcal{B})$ be the subalgebra of $\mathcal{H}(\mathcal{D})$ generated by $u_E,E\in\mathcal{B}$. Since $\mathcal{B}$ is extension closed, we have
\[\mathcal{H}(\mathcal{B})=\bigoplus_{\gamma\in\Lambda_\phi}\mathcal{H}(\mathcal{B})_\gamma,\]
with each $\mathcal{H}(\mathcal{B})_\gamma$ being finite dimensional, by Lemma \ref{sector degree finite}. The following lemma is a consequence of Lemma \ref{coincide Ringel-Hall}.
\begin{lemma}
    $\mathcal{H}(\mathcal{B})$ coincides with the Ringel-Hall algebra of $\mathcal{B}$.
\end{lemma}

Let
\[\hat{\mathcal{H}}_\sigma(\phi)=\prod_{\gamma\in\Lambda_\phi}\mathcal{H}(\mathcal{B})_\gamma\]
be the completion of $\mathcal{H}(\mathcal{B})$, in the sense of Proposition \ref{sector completion}. Consider the positively graded part
\[\mathfrak{g}_\phi=\prod_{\gamma\in\Lambda_\phi^+}\mathcal{H}(\mathcal{B})_\gamma.\]
For any $R>0$, we set 
\[\Lambda_{\phi,R}=h_\phi^{-1}((0,R])\subseteq\Lambda_\phi^+,\]
and
\[F^{>R}\mathfrak{g}_\phi=\prod_{\gamma\in\Lambda_\phi^+-\Lambda_{\phi,R}}\mathcal{H}(\mathcal{B})_\gamma.\]
It is an ideal of $\mathfrak{g}_\phi$ and the corresponding quotient algebra is
\[\mathfrak{g}_\phi/F^{>R}\mathfrak{g}_\phi=\bigoplus_{\gamma\in\Lambda_{\phi,R}}\mathcal{H}(\mathcal{B})_\gamma,\]
which is a nilpotent algebra. When $R_1>R_2>0$, we have the canonical inclusion and projection
\[i_{R_1,R_2}:F^{>R_1}\mathfrak{g}_\phi\hookrightarrow F^{>R_2}\mathfrak{g}_\phi,p_{R_1,R_2}:\mathfrak{g}_\phi/F^{>R_1}\mathfrak{g}_\phi\twoheadrightarrow\mathfrak{g}_\phi/F^{>R_2}\mathfrak{g}_\phi.\]
It is easy to check that
\[\varprojlim_{R>0}\mathfrak{g}_\phi/F^{>R}\mathfrak{g}_\phi=\mathfrak{g}_\phi.\]
Now we define
\begin{align*}
    \exp_R:\mathfrak{g}_\phi/F^{>R}\mathfrak{g}_\phi&\rightarrow u_0+\mathfrak{g}_\phi/F^{>R}\mathfrak{g}_\phi\\
        x&\mapsto\sum_{n=0}^\infty\dfrac{1}{n!}x^n.
\end{align*}
It is well-defined, by the nilpotency of $\mathfrak{g}_\phi/F^{>R}\mathfrak{g}_\phi$. One can easily check that the diagram
\begin{center}
    \begin{tikzcd}
\mathfrak{g}_\phi/F^{>R_1}\mathfrak{g}_\phi \arrow[r, "\exp_{R_1}"] \arrow[d, "{p_{R_1,R_2}}"] & u_0+\mathfrak{g}_\phi/F^{>R_1}\mathfrak{g}_\phi \arrow[d, "{p_{R_1,R_2}}"] \\
\mathfrak{g}_\phi/F^{>R_2}\mathfrak{g}_\phi \arrow[r, "\exp_{R_2}"]                            & u_0+\mathfrak{g}_\phi/F^{>R_2}\mathfrak{g}_\phi                           
\end{tikzcd}
\end{center}
commutes for any $R_1>R_2>0$. So we obtain the \emph{exponential map}
\begin{align*}
    \exp:\mathfrak{g}_\phi&\rightarrow u_0+\mathfrak{g}_\phi\\
    x&\mapsto\sum_{n=0}^\infty\dfrac{1}{n!}x^n.
\end{align*}
Similarly, if we set
\begin{align*}
    \log_R:u_0+\mathfrak{g}_\phi/F^{>R}\mathfrak{g}_\phi&\rightarrow \mathfrak{g}_\phi/F^{>R}\mathfrak{g}_\phi\\
        y&\mapsto\sum_{n=1}^\infty\dfrac{(-1)^{n-1}}{n}(y-u_0)^n,
\end{align*}
then it can be lifted to the \emph{logarithm map}
\begin{align*}
    \log:u_0+\mathfrak{g}_\phi&\rightarrow \mathfrak{g}_\phi\\
        y&\mapsto\sum_{n=1}^\infty\dfrac{(-1)^{n-1}}{n}(y-u_0)^n.
\end{align*}
By the relationship between the exponential and logarithm maps, we have
\begin{proposition}
    The maps $\exp_R$ and $\log_R$ are inverses of each other, and both are bijections between $\mathfrak{g}_\phi/F^{>R}\mathfrak{g}_\phi$ and $u_0+\mathfrak{g}_\phi/F^{>R}\mathfrak{g}_\phi$ for any $R>0$. As a consequence, $\exp$ and $\log$ are inverses of each other between $\mathfrak{g}_\phi$ and $u_0+\mathfrak{g}_\phi$.
\end{proposition}
\begin{corollary}
    The set
    \[G_\phi=u_0+\mathfrak{g}_\phi\subseteq\hat{\mathcal{H}}_\sigma(\phi)^\times\]
    is a pro-nilpotent group.
\end{corollary}
\begin{proof}
    Indeed,
    \begin{align*}
        G_\phi=u_0+\mathfrak{g}_\phi&=u_0+\varprojlim_{R>0}\mathfrak{g}_\phi/F^{>R}\mathfrak{g}_\phi\\
        &=\varprojlim_{R>0}(u_0+\mathfrak{g}_\phi/F^{>R}\mathfrak{g}_\phi)\\
        &=\varprojlim_{R>0}\exp_R(\mathfrak{g}_\phi/F^{>R}\mathfrak{g}_\phi).
    \end{align*}
\end{proof}
\subsection{Restrictions of nearby stability conditions}
Let $R>0$. Choose $0<\eta_R<\dfrac{1}{4}$ as in Proposition \ref{nearby}. We fix a neighbourhood of $\sigma$ in $\Stab_\cl(\mathcal{D})$:
\[V_R(\sigma)=\{\tau\in\Stab_\cl(\mathcal{D})|d(\sigma,\tau)<\eta_R\},\]
whose elements will be called \emph{nearby stability conditions} or \emph{near} $\sigma$ with respect to $R$. Let $\tau$ be near $\sigma$ with respect to some $R>0$. The function
\[Z_\tau|_\mathcal{B}:\Lambda_\phi\rightarrow\mathbb{C}\]
is additive and $Z_\tau(\Lambda_\phi^+)$ lies in $e^{i\pi(\phi-\frac12)}H$. Thus $Z_\tau|_\mathcal{B}$ is a central charge on $\Lambda_\phi$. By the results in the last subsection and Lemma \ref{ordered decomposition} we obtain a bijection
\begin{align*}
    \Phi_\tau:\mathfrak{g}_\phi&\rightarrow G_\phi\\
    x=(x_\gamma)&\mapsto\prod_\psi^\leftarrow\exp\left(\sum_{\phi_\tau(\gamma)=\psi}x_\gamma\right).
\end{align*}
Our study will focus on the element
\[\Phi_\tau^{-1}(\1_{\mathcal{B}})\in\mathfrak{g}_\phi.\]
\begin{lemma}\label{semistable agree}
    $Z_\tau|_{\mathcal{B}}$ is a rotation of a stability function on $\mathcal{B}$ with the HN property. Moreover, if $0\neq E\in\mathcal{B}$ satisfies $m_\sigma(E)=h_\phi(\cl E)\leqslant R$, then it is $Z_\tau|_{\mathcal{B}}$-semistable if and only if it is $\tau$-semistable.
\end{lemma}
\begin{proof}
    The discussion above has shown that $Z_\tau|_{\mathcal{B}}$ is a rotation of a stability function on $\mathcal{B}$. Since $\mathcal{B}$ is of finite length, it must have the HN property.

    Let $0\neq E\in\mathcal{B}$ with $m_\sigma(E)\leqslant R$. By Proposition \ref{nearby}(2), each $\tau$-HN factor of $E$ belongs to $\mathcal{B}$, and they thus form an HN filtration of $E$ in $\mathcal{B}$ under $Z_\tau|_{\mathcal{B}}$. By including $\mathcal{B}$ in a rotated heart $\mathcal{P}_\tau\left(\phi-\dfrac12,\phi+\dfrac12\right]$ and using Lemma \ref{heart}, we obtain the result.
\end{proof}
\begin{corollary}\label{delta agree}
    Let $\delta^\tau_\gamma\in\mathcal{H}(\mathcal{B})\subseteq\mathcal{H}(\mathcal{D})$ be the characteristic function supported on $Z_\tau|_{\mathcal{B}}$-semistable objects in $\mathcal{B}$ with $\cl E=\gamma\in\Lambda_{\phi,R}$. Then
    \[\delta^\tau_\gamma=\delta^\tau_{\phi_\tau(\gamma),\gamma}\in\mathcal{H}(\mathcal{D}).\]
\end{corollary}
\begin{proof}
    By Lemma \ref{semistable agree}, we only need to show that if $0\neq E\in\mathcal{P}_\tau(\phi_\tau(\gamma))$ with $\cl E=\gamma\in\Lambda_{\phi,R}$, then $E\in\mathcal{B}$.

    To show that $E$ is $\sigma$-semistable, by Proposition \ref{nearby}, we only need to show that $m_\tau(E)\leqslant2R$. Let $c$ be a support constant of $\sigma$. By the proof of Corollary \ref{finite distance}, $\tau$ has a support constant $e^{-d(\sigma,\tau)}c$. By the choice of $\eta_R$, we have
    \[e^{-d(\sigma,\tau)}c\geqslant\dfrac{c}{2},\|Z_\sigma-Z_\tau\|\leqslant\dfrac{c}{4}.\]
    Since $E$ is $\tau$-semistable, the support property of $\tau$ implies
    \[\|\gamma\|\leqslant\dfrac{2}{c}|Z_\tau(\gamma)|.\]
    By the triangle inequality,
    \begin{align*}
        R\geqslant h_\phi(\gamma)&=|Z_\sigma(\gamma)|\geqslant|Z_\tau(\gamma)|-\|Z_\tau-Z_\sigma\|\|\gamma\|\geqslant|Z_\tau(\gamma)|-\dfrac{c}{4}\dfrac{2}{c}|Z_\tau(\gamma)|\\
        &=\dfrac12|Z_\tau(\gamma)|=\dfrac12m_\tau(E),
    \end{align*}
    that is, $m_\tau(E)\leqslant2R$.

    Now $E$ is $\sigma$-semistable with $\cl E\in\Lambda_{\phi,R}$, so we must have
    \[\phi_\sigma(E)=\phi+2k\]
    for some $k\in\mathbb{Z}$. However, $d(\sigma,\tau)<\eta_R$ forces that $k=0$. So $E\in\mathcal{P}_\sigma(\phi)=\mathcal{B}$.
\end{proof}
As a consequence, if $\gamma\in\Lambda_{\phi,R}$, then the Joyce element
\[\epsilon^\tau_\gamma=\sum_{\begin{subarray}{c}
    \gamma_1+\dots+\gamma_n=\gamma,\gamma_i\in\Lambda_\phi^+,n\geqslant1\\
    \phi_\tau(\gamma_i)=\phi_\tau(\gamma)
\end{subarray}}\dfrac{(-1)^{n-1}}{n}\delta^\tau_{\gamma_1}*\dots*\delta^\tau_{\gamma_n}\]
with respect to the stability function $Z_\tau|_{\mathcal{B}}$, it equals $\epsilon^\tau_{\phi_\tau(\gamma),\gamma}$ as defined in (\ref{Joyce's element}), since each $\tau$-semistable object $E\in\mathcal{P}_\tau(\phi_\tau(\gamma))$ with $\cl E=\gamma_i$ must belong to $\mathcal{B}$, by Proposition \ref{nearby}. Then we have:
\begin{proposition}\label{decomposition of 1}
    We write $\1_R$ for the projection of $\1_{\mathcal{B}}$ in $u_0+\mathfrak{g}_\phi/F^{>R}\mathfrak{g}_\phi$. Then for any $\tau$ near $\sigma$ with respect to $R$, we have
    \[\1_R=\prod_\psi^\leftarrow\exp_R\left(\sum_{\gamma\in\Lambda_{\phi,R},\phi_\tau(\gamma)=\psi}\epsilon^\tau_\gamma\right).\]
\end{proposition}
\begin{proof}
    Let us first show that
    \[\1_R=\prod_{\psi\in(\phi-\eta_R,\phi+\eta_R)}^\leftarrow\left(u_0+\sum_{\gamma\in\Lambda_{\phi,R},\phi_\tau(
    \gamma)=\psi}\delta^\tau_\gamma\right)\]
    in $u_0+\mathfrak{g}_\phi/F^{>R}\mathfrak{g}_\phi$.
    We show it by expanding the right-hand side:
    \begin{align*}
    \RHS&=\prod_{\psi\in(\phi-\eta_R,\phi+\eta_R)}^\leftarrow\left(u_0+\sum_{\gamma\in\Lambda_{\phi,R},\phi_\tau(
    \gamma)=\psi}\delta^\tau_\gamma\right)\\
    &=u_0+\sum_{\begin{subarray}{c}
        \gamma_1,\dots,\gamma_r\in\Lambda_{\phi,R}\\
        \phi+\eta_R>\phi_\tau(\gamma_1)>\dots>\phi_\tau(\gamma_r)>\phi-\eta_R\\
    \end{subarray}} \delta^\tau_{\gamma_1}*\dots*\delta^\tau_{\gamma_r}\\
    &=u_0+\sum_{\gamma\in\Lambda_\phi^+}\sum_{\begin{subarray}{c}
        \gamma_1,\dots,\gamma_r\in\Lambda_{\phi,R},\gamma_1+\dots+\gamma_r=\gamma\\
        \phi+\eta_R>\phi_\tau(\gamma_1)>\dots>\phi_\tau(\gamma_r)>\phi-\eta_R\\
    \end{subarray}} \delta^\tau_{\gamma_1}*\dots*\delta^\tau_{\gamma_r}\\
    &\equiv u_0+\sum_{\gamma\in\Lambda_\phi^+}\sum_{\begin{subarray}{c}
        \gamma_1+\dots+\gamma_r=\gamma\\
        \phi+\eta_R>\phi_\tau(\gamma_1)>\dots>\phi_\tau(\gamma_r)>\phi-\eta_R\\
    \end{subarray}} \delta^\tau_{\gamma_1}*\dots*\delta^\tau_{\gamma_r}\,(\mmod F^{>R}\mathfrak{g}_\phi)\\
    (\ref{HN identity in sector})&=u_0+\sum_{\gamma\in\Lambda_\phi^+}\kappa^\tau_{(\phi-\eta_R,\phi+\eta_R),\gamma}\\
    &=\1_\mathcal{B}.
    \end{align*}
    That is,
    \[\1_\mathcal{B}\equiv\prod_{\psi\in(\phi-\eta_R,\phi+\eta_R)}^\leftarrow\left(u_0+\sum_{\gamma\in\Lambda_{\phi,R},\phi_\tau(
    \gamma)=\psi}\delta^\tau_\gamma\right)\mod F^{>R}\mathfrak{g}_\phi\]
    in $u_0+\mathfrak{g}_\phi$, which is the same as we want. The rest is to show that
    \[u_0+\sum_{\gamma\in\Lambda_{\phi,R},\phi_\tau(
    \gamma)=\psi}\delta^\tau_\gamma=\exp_R\left(\sum_{\gamma\in\Lambda_{\phi,R},\phi_\tau(\gamma)=\psi}\epsilon^\tau_\gamma\right)\]
    in $u_0+\mathfrak{g}_\phi/F^{>R}\mathfrak{g}_\phi$ for any $\psi\in(\phi-\eta_R,\phi+\eta_R)$. By taking the logarithm of the left hand side, we obtain
    \begin{align*}
        \log\LHS&=\log\left(u_0+\sum_{\gamma\in\Lambda_{\phi,R},\phi_\tau(
    \gamma)=\psi}\delta^\tau_\gamma\right)\\
    &=\sum_{n=1}^\infty\dfrac{(-1)^{n-1}}{n}\left(\sum_{\gamma\in\Lambda_{\phi,R},\phi_\tau(
    \gamma)=\psi}\delta^\tau_\gamma\right)^n\\
    &=\sum_{n=1}^\infty\dfrac{(-1)^{n-1}}{n}\sum_{\begin{subarray}{c}
        \gamma_1,\dots,\gamma_n\in\Lambda_{\phi,R}\\
        \phi_\tau(\gamma_1)=\dots=\phi_\tau(\gamma_n)=\psi
    \end{subarray}}\delta^\tau_{\gamma_1}*\dots*\delta^\tau_{\gamma_n}\\
    &=\sum_{\begin{subarray}{c}
        \gamma_1,\dots,\gamma_n\in\Lambda_{\phi,R},n\geqslant1\\
        \phi_\tau(\gamma_1)=\dots=\phi_\tau(\gamma_n)=\psi
    \end{subarray}}\dfrac{(-1)^{n-1}}{n}\delta^\tau_{\gamma_1}*\dots*\delta^\tau_{\gamma_n}\\
    &=\sum_{\gamma\in\Lambda_\phi^+}\sum_{\begin{subarray}{c}
        \gamma_1,\dots,\gamma_n\in\Lambda_{\phi,R}\\
        \gamma_1+\dots+\gamma_n=\gamma,n\geqslant1\\
        \phi_\tau(\gamma_1)=\dots=\phi_\tau(\gamma_n)=\psi
    \end{subarray}}\dfrac{(-1)^{n-1}}{n}\delta^\tau_{\gamma_1}*\dots*\delta^\tau_{\gamma_n}\\
    &\equiv\sum_{\gamma\in\Lambda_{\phi,R}}\sum_{\begin{subarray}{c}
        \gamma_1+\dots+\gamma_n=\gamma,n\geqslant1\\
        \phi_\tau(\gamma_1)=\dots=\phi_\tau(\gamma_n)=\psi
    \end{subarray}}\dfrac{(-1)^{n-1}}{n}\delta^\tau_{\gamma_1}*\dots*\delta^\tau_{\gamma_n}\,(\mmod F^{>R}\mathfrak{g}_\phi)\\
    (\ref{Joyce's element})&=\sum_{\gamma\in\Lambda_{\phi,R},\phi_\tau(\gamma)=\psi}\epsilon^\tau_\gamma.
    \end{align*}
    We complete the proof by applying the exponential map on both sides.
\end{proof}
Let $\tau_+,\tau_-$ be two stability conditions near $\sigma$ with respect to $R$. Proposition \ref{decomposition of 1} implies
\[\prod_\psi^\leftarrow\exp_R\left(\sum_{\gamma\in\Lambda_{\phi,R},\phi_{\tau_+}(\gamma)=\psi}\epsilon^{\tau_+}_\gamma\right)=\1_R=\prod_\psi^\leftarrow\exp_R\left(\sum_{\gamma\in\Lambda_{\phi,R},\phi_{\tau_-}(\gamma)=\psi}\epsilon^{\tau_-}_\gamma\right).\]
By applying \cite[Theorem 5.2]{Joy2008} and Theorem \ref{universal Lie}, we obtain the following result:
\begin{theorem}\label{local Joyce transformation}
    For any $\gamma\in\Lambda_{\phi,R}$, we have an identity
    \begin{equation}\label{Joyce transformation}
        \epsilon^{\tau_+}_\gamma=\sum_{\gamma_1+\dots+\gamma_n=\gamma,n\geqslant1}U(\gamma_1,\dots,\gamma_n;\tau_-,\tau_+)\epsilon^{\tau_-}_{\gamma_1}*\dots*\epsilon^{\tau_-}_{\gamma_n},
    \end{equation}
    where $U(\gamma_1,\dots,\gamma_n;\tau_-,\tau_+)\in\mathbb{Q}$ is the Joyce transformation coefficient on $\epsilon$. Moreover, the right-hand side of (\ref{Joyce transformation}) is in fact a Lie polynomial in $\epsilon^{\tau_-}_{\gamma_1},\dots,\epsilon^{\tau_-}_{\gamma_n}$.
\end{theorem}
\subsection{Finite numerical Joyce Transport}
Let $\Sigma$ be a connected component of $\Stab_\cl(\mathcal{D})$. Suppose there exists one $\sigma\in\Sigma$ satisfying the support property and the finite counting condition.

\begin{lemma}\label{one point to the component}
    Any $\tau\in\Sigma$ satisfies the support property and the finite counting condition.
\end{lemma}
\begin{proof}
    Since $\Sigma$ is connected in the complex manifold $\Stab_\cl(\mathcal{D})$, it must be path connected. Pick a path $\alpha:[0,1]\rightarrow\Sigma$ with $\alpha_0=\sigma$ and $\alpha_1=\tau$. The compactness of $\alpha([0,1])$ and the continuity of the distance function $d$ imply that $d(\sigma,\tau)<+\infty$. Hence, the result follows by Corollary \ref{finite distance}.
\end{proof}
Let $\sigma\in\Sigma$. We set $\mathfrak{L}_\sigma$ to be the Lie subalgebra of $\mathcal{H}(\mathcal{D})$ generated by the subset
\[\{\epsilon^\sigma_{\phi,\gamma}|\phi\in\mathbb{R},0\neq\gamma\in\Gamma\}.\]
Now we can state the main theorem of this article.
\begin{theorem}[Joyce transport]\label{Joyce transport}
    Let $\sigma_0,\sigma_1\in\Sigma$. For any $\phi\in\mathbb{R},\gamma\in\Gamma$, the Joyce element $\epsilon^{\sigma_1}_{\phi,\gamma}$ can be written as a finite Lie polynomial in $\epsilon^{\sigma_0}_{\psi,\beta},\psi\in\mathbb{R},\beta\in\Gamma$ with coefficients in $\mathbb{Q}$, with each Lie monomial having total degree $\gamma$. As a consequence, we have
    \[\mathfrak{L}_{\sigma_0}=\mathfrak{L}_{\sigma_1}.\]
\end{theorem}
\begin{proof}
    The spirit of the proof is to pick a path from $\sigma_0$ to $\sigma_1$ and transform all the Joyce elements along the path. Without loss of generality, we may assume that $\epsilon^{\sigma_1}_{\phi,\gamma}\neq0$. That forces that $Z_{\sigma_1}(\gamma)$ has phase $\phi+2k$ for some $k\in\mathbb{Z}$. We divide the proof into steps.
    
    \emph{Step 1: Choose a path and a height cutoff.} Let $\sigma:[0,1]\rightarrow\Sigma$ be a $C^1$ path in $\Sigma$ from $\sigma_0$ to $\sigma_1$. For each $t\in[0,1]$, we write $\sigma_t=(Z_t,\mathcal{P}_t)$. We denote by $\|\cdot\|$ the operator norm on $\Gamma_{\mathbb{C}}^*$ induced by the norm on $\Gamma_{\mathbb{C}}$. Let $c_0$ be a support constant of $\sigma_0$. By the compactness of $[0,1]$ and the continuity of $t\mapsto d(\sigma_0,\sigma_t)$, $M=\displaystyle{\max_{t\in[0,1]}d(\sigma_0,\sigma_t)}$ exists. By Corollary \ref{finite distance}, $c=e^{-M}c_0$ is a common support constant of $\sigma_t$ for any $t\in[0,1]$. Now set
    \[V=\int_0^1\left\|\dfrac{\diff Z_t}{\diff t}\right\|\diff t,B=e^{\frac{4V}{c}}|Z_1(\gamma)|,\]
    and fix $R=2B$.

    \emph{Step 2: Divide $[0,1]$ into finite segments.} For any $t\in[0,1]$, we can pick an open segment $J_t=(c_t,d_t)\subseteq\mathbb{R}$, such that:
    \begin{itemize}
        \item For any $s\in J_t\cap[0,1]$, we have $\sigma_s\in V_R(\sigma_t)$. It is worth warning that the radius of each $V_R(\sigma_t)$ might be different.
        \item The variation of the central charge is small enough:
        \[\int_{J_t\cap[0,1]}\left\|\dfrac{\diff Z_s}{\diff s}\right\|\diff s<\dfrac{c}{2}.\]
    \end{itemize}
    By the compactness of $[0,1]$, one can choose a finite subset $\{t_1,\dots,t_N\}\subseteq[0,1]$, such that
    \[[0,1]\subseteq\bigcup_{j=1}^N J_{t_j}=\bigcup_{j=1}^N(c_{t_j},d_{t_j}).\]
    We can shrink each segment $(c_{t_j},d_{t_j})$ such that $t_j<c_{t_{j+1}}<d_{t_j}<t_{j+1},j=1,\dots,N-1$. Then we can pick $0=a_0<a_1<\dots<a_N=1$, such that $t_j\in[a_{j-1},a_j]\subseteq J_{t_j}$ for each $j=1,\dots,N$. Now we have
    \[x_j:=\dfrac{1}{c}\int_{a_{j-1}}^{a_j}\left\|\dfrac{\diff Z_t}{\diff t}\right\|\diff t<\dfrac{1}{2},j=1,\dots,N.\]

    \emph{Step 3: Calculate on each segment.} Now we focus on a small segment $[a_{j-1},a_j]$. By the hypothesis, $\sigma_{a_{j-1}},\sigma_{a_j}\in V_R(\sigma_{t_j})$. Let $\psi\in\mathbb{R}$, and write
    \[\Lambda^j_{\psi,+}=\{\cl E|0\neq E\in\mathcal{P}_{t_j}(\psi)\},\Lambda^j_{\psi,+}=\{\cl E|0\neq E\in\mathcal{P}_{t_j}(\psi),m_{t_j}(E)\leqslant R\}.\]
    For any $\beta\in\Lambda^j_{\psi,+}$ and $s\in[a_{j-1},a_j]$, we have
    \[|Z_s(\beta)-Z_{t_j}(\beta)|\leqslant\|Z_s-Z_{t_j}\|\|\beta\|\leqslant cx_j\|\beta\|\leqslant x_j|Z_{t_j}(\beta)|.\]
    So
    \begin{equation}\label{1-x}
        (1-x_j)|Z_{t_j}(\beta)|\leqslant|Z_s(\beta)|\leqslant(1+x_j)|Z_{t_j}(\beta)|.
    \end{equation}
    Now consider the Joyce element $\epsilon^{\sigma_{a_j}}_\alpha=\epsilon^{\sigma_{a_j}}_{\theta,\alpha}$, where $\alpha=\cl E$ for some $0\neq E\in\mathcal{P}_{a_j}(\theta)$ with $|Z_{a_j}(\alpha)|\leqslant B$. $\epsilon^{\sigma_{a_j}}_\alpha$ is a function supported on $\sigma_{a_j}$-semistable objects with $\sigma_{a_j}$-mass no more than $B$. By Proposition \ref{nearby}, they are all $\sigma_{t_j}$-semistable. By the inequality (\ref{1-x}), their $\sigma_{t_j}$-masses are no more than $R=2B$. Now we apply Theorem \ref{local Joyce transformation} to $\epsilon^{\sigma_{a_j}}_\alpha$, so $\epsilon^{\sigma_{a_j}}_\alpha$ can be written as a Lie polynomial in several $\epsilon^{\sigma_{a_{j-1}}}_{\beta}$'s, where for each Lie monomial in $\epsilon^{\sigma_{a_{j-1}}}_{\beta_1},\dots,\epsilon^{\sigma_{a_{j-1}}}_{\beta_r}$, we have $\beta_1+\dots+\beta_r=\alpha$ and their central charges under $Z_{t_j}$ lie on the same positive ray. Thus we have
    \[|Z_{{t_j}}(\alpha)|=\left|\sum_{i=1}^r Z_{{t_j}}(\beta_i)\right|=\sum_{i=1}^r|Z_{{t_j}}(\beta_i)|.\]
    When $0<x<\dfrac{1}{2}$, one can easily show that
    \[\log\dfrac{1+x}{1-x}\leqslant 4x.\]
    By the inequality (\ref{1-x}), we have
    \[\sum_{i=1}^r|Z_{{a_{j-1}}}(\beta_i)|\leqslant\dfrac{1+x_j}{1-x_j}|Z_{{a_j}}(\alpha)|\leqslant e^{4x_j}|Z_{{a_j}}(\alpha)|.\]
    In conclusion, for any $\beta$ such that the Joyce element $\epsilon^{\sigma_{a_{j-1}}}_\beta$ occurs in the Lie polynomial expression of $\epsilon^{\sigma_{a_j}}_\alpha$, we must have
    \begin{equation}\label{e^4x}
        |Z_{{a_{j-1}}}(\beta)|\leqslant e^{4x_j}|Z_{{a_j}}(\alpha)|.
    \end{equation}

    \emph{Step 4: Express backwards along the path.} We find the expression of $\epsilon^{\sigma_1}_\gamma$ inductively backwards along the path $\sigma$. By (\ref{e^4x}), each $\beta$ such that $\epsilon^{\sigma_{a_{N-1}}}_\beta$ occurs in the Lie polynomial expression of $\epsilon^{\sigma_1}_\gamma$ satisfies
    \[|Z_{a_{N-1}}(\beta)|\leqslant e^{4x_N}|Z_1(\gamma)|.\]
    Inductively, each $\beta$ such that $\epsilon^{\sigma_{a_{j-1}}}_\beta$ occurs in the Lie polynomial expression of $\epsilon^{\sigma_1}_\gamma$ must satisfy
    \begin{align*}
        |Z_{a_{j-1}}(\beta)|&\leqslant e^{4(x_j+\dots+x_N)}|Z_1(\gamma)|\\
        &=\exp\left(\dfrac{4}{c}\int_{a_{j-1}}^1\left\|\dfrac{\diff Z_t}{\diff t}\right\|\diff t\right)|Z_1(\gamma)|\\
        &\leqslant\exp\left(\dfrac{4V}{c}\right)|Z_1(\gamma)|\\
        &=B.
    \end{align*}
    Thus the expression process can inductively go back to $\sigma_0$ and we finish the proof.
\end{proof}
\section{Derived Hall scattering groupoids}
In this section, we use Bridgeland's method of scattering diagrams to realize the process of Joyce transport. For the details of the notion and basic properties of Bridgeland scattering diagrams, we refer the reader to \cite{Bri2017}. We only provide a brief recall for convenience.
\subsection{Bridgeland scattering diagram on phase categories}
Let $\sigma\in\Stab_\cl(\mathcal{D})$ satisfy the support property and the finite counting condition. We keep the notation as before. Let $\Gamma_{\mathcal{B}}$ be the subgroup of $\Gamma$ generated by $\Lambda_\phi$, and set
\[M_\mathcal{B}=\Hom_{\mathbb{Z}}(\Gamma_\mathcal{B},\mathbb{R}).\]
One can check that the height function $h_\phi$ can be extended to a function in $M_\mathcal{B}$, and we keep the notation $h_\phi$ for the extension. The proof of the following lemma is similar to that of Lemma \ref{semistable agree}.
\begin{lemma}
    For any $\theta\in M_\mathcal{B}$, the function
    \begin{align*}
        Z_\theta:\Lambda_\phi^+&\rightarrow\mathbb{C}\\
        \beta&\mapsto-\theta(\beta)+ih_\phi(\beta)
    \end{align*}
    defines a stability function on $\mathcal{B}$ with HN property.
\end{lemma}

Given $\theta\in M_\mathcal{B}$, we set the following three subcategories of $\mathcal{B}$:
\begin{align*}
    \mathcal{C}_\theta&=\{E\in\mathcal{B}|\theta(\cl E)=0,\theta(\cl A)\leqslant0\text{ for any }A\subseteq E\},\\
    \mathcal{T}_\theta&=\{E\in\mathcal{B}|\theta(\cl Q)>0\text{ for any nonzero quotient }Q\text{ of }E\},\\
   \overline{ \mathcal{T}}_\theta&=\{E\in\mathcal{B}|\theta(\cl Q)\geqslant0\text{ for any nonzero quotient }Q\text{ of }E\}.
\end{align*}
The following lemma follows immediately from the definition:
\begin{lemma}\label{T perp C}
    If $T\in\mathcal{T}_\theta$ and $C\in\mathcal{C}_\theta$, then $\Hom_\mathcal{B}(T,C)=0$.
\end{lemma}
\begin{proof}
    Let $f\in\Hom_\mathcal{B}(T,C)$. If $Q=\im f\neq 0$, since $T\in\mathcal{T}_\theta$, we must have $\theta(\cl Q)>0$. However, $Q$ is a subobject of $C\in\mathcal{C}_\theta$, so $\theta(\cl Q)\leqslant0$, a contradiction.
\end{proof}
We define the \emph{slope function} corresponding to $\theta$ to be
\begin{align*}
    \mu_\theta:\mathcal{B}-\{0\}&\rightarrow\mathbb{R}\\
    E&\mapsto\dfrac{\theta(\cl E)}{h_\phi(\cl E)}.
\end{align*}
Then the three subcategories have the following description:
\begin{lemma}\label{torsion part}
    Let $\theta\in M_\mathcal{B}$ and $0\neq E\in\mathcal{B}$.
    \begin{enumerate}
        \item $E\in\overline{\mathcal{T}}_\theta$ if and only if the slope of each $Z_\theta$-HN factor of $E$ is nonnegative.
        \item $E\in\mathcal{T}_\theta$ if and only if the slope of each $Z_\theta$-HN factor of $E$ is positive.
    \end{enumerate}
\end{lemma}
\begin{proof}
    We only prove (1), since the proof of (2) is similar. Recall that if there is a short exact sequence $0\rightarrow A\rightarrow F\rightarrow Q\rightarrow0$ with $A,F,Q\neq0$ and $F$ being $Z_\theta$-semistable, then we have
    \[\mu_\theta(A)\leqslant\mu_\theta(F)\leqslant\mu_\theta(Q).\]
    Let $0=E_0\subsetneqq E_1\subsetneqq\dots\subsetneqq E_m=E$ be the HN filtration of $E$ under $Z_\theta$. The HN factors are $F_i=E_i/E_{i-1},i=1,\dots,m$, with $\mu_\theta(F_1)>\dots>\mu_\theta(F_m)$.

    For the only-if part, note that $F_m=E/E_{m-1}$ is a nonzero quotient of $E$. Since $E\in\overline{\mathcal{T}}_\theta$, we have $\mu_\theta(F_m)\geqslant0$, and hence
    \[\mu_\theta(F_1)>\dots>\mu_\theta(F_m)\geqslant0.\]

    For the if part, let $Q$ be a nonzero quotient of $E$ with the canonical projection $\pi:E\twoheadrightarrow Q$. Set $Q_i=\pi(E_i)$. Then $0=Q_0\subseteq Q_1\subseteq\dots\subseteq Q_m=Q$ is a filtration of $Q$. The quotient factor $G_i=Q_i/Q_{i-1}$ is a quotient of $F_i,i=1,\dots,m$. If $G_i\neq0$, since $F_i$ is $Z_\theta$-semistable, we have $\mu_\theta(G_i)\geqslant\mu_\theta(F_i)\geqslant0$, and therefore $\theta(\cl G_i)\geqslant0$. If $G_i=0$, then $\theta(\cl G_i)=0$. So
    \[\theta(\cl Q)=\sum_{i=1}^m\theta(\cl G_i)\geqslant0,\]
    and hence $\mu_\theta(Q)\geqslant0$.
\end{proof}

\begin{lemma}\label{torsion semistable decomposition}
    Let $\theta\in M_\mathcal{B}$ and $E\in\mathcal{B}$. $E\in\overline{\mathcal{T}}_\theta$ if and only if there exists a short exact sequence
    \[0\rightarrow T\rightarrow E\rightarrow C\rightarrow0\]
    with $T\in\mathcal{T}_\theta$ and $C\in\mathcal{C}_\theta$. Moreover, the sequence is unique up to unique isomorphism.
\end{lemma}
\begin{proof}
    When $E\in\mathcal{T}_\theta$ or $E\in\mathcal{C}_\theta$, there is nothing to prove. So let us assume that $E\notin\mathcal{T}_\theta$ and $E\notin\mathcal{C}_\theta$. Let $0=E_0\subsetneqq E_1\subsetneqq\dots\subsetneqq E_m=E$ be the HN filtration of $E$ under $Z_\theta$. The HN factors are $F_i=E_i/E_{i-1},i=1,\dots,m$, with $\mu_\theta(F_1)>\dots>\mu_\theta(F_m)$.

    If $E\in\overline{\mathcal{T}}_\theta$, by Lemma \ref{torsion part}(1), we have $\mu_\theta(F_1)>\dots>\mu_\theta(F_m)\geqslant0$. If $\mu_\theta(F_m)=0$, pick $T=E_{m-1}$ and $C=F_m$. If $\mu_\theta(F_m)>0$, pick $T=E$ and $C=0$.

    Conversely, suppose that such a short exact sequence exists. The cases when $T=0$ or $C=0$ are easy, so let us only show the case that $T\neq0,C\neq0$. Let $0=T_0\subsetneqq T_1\subsetneqq\dots\subsetneqq T_k=T$ be the HN filtration of $T$. By Lemma \ref{torsion part}(2), we have $\mu_\theta(T_1/T_0)>\dots>\mu_\theta(T_k/T_{k-1})>0=\mu_\theta(C)=\mu_\theta(E/T)$. Thus
    \[0=T_0\subsetneqq T_1\subsetneqq\dots\subsetneqq T_k=T\subsetneqq E\]
    is the HN filtration of $E$. Then $E\in\overline{\mathcal{T}}_\theta$ by Lemma \ref{torsion part}(1).

    For the uniqueness part, let
    \begin{align*}
    0\rightarrow T\xrightarrow{i} &E\xrightarrow{p} C\rightarrow0\\
        0\rightarrow T'\xrightarrow{i'} &E\xrightarrow{p'} C'\rightarrow0
    \end{align*}
     be two short exact sequences in $\mathcal{B}$ with $T,T'\in\mathcal{T}_\theta,C,C'\in\mathcal{C}_\theta$. Applying the left exact contravariant functor $\Hom_\mathcal{B}(-,C)$ to the second sequence, we obtain the exact sequence
    \[0\rightarrow\Hom(C',C)\xrightarrow{-\circ p'}\Hom(E,C)\rightarrow\Hom(T',C).\]
    Since $\Hom_\mathcal{B}(T',C)=0$ by Lemma \ref{T perp C}, $-\circ p'$ is an isomorphism. Thus, there exists a unique morphism $g':C'\rightarrow C$, such that $g'\circ p'=p$. We omit the rest of the argument since it is quite similar to that of Proposition \ref{phase truncation}.
\end{proof}
We translate Lemma \ref{torsion semistable decomposition} into the language of Hall algebras:
\begin{corollary}\label{torsion semistable decomposition, hall}
    In the pro-nilpotent group $G_\phi$, we have the identity
    \[\1_{\overline{\mathcal{T}}_\theta}=\1_{\mathcal{T}_\theta}*\1_{\mathcal{C}_\theta}.\]
    The identity also holds in any cutoff $G_{\phi,R}=u_0+\mathfrak{g}_\phi/F^{>R}\mathfrak{g}_\phi$.
\end{corollary}
We set
\[\Xi_R=\{\cl F| F\neq 0\text{ is a subquotient of some $E$ in $\mathcal{B}$ such that }h_\phi(\cl E)\leqslant R\}.\]
It is a subset of $\Lambda_{\phi,R}$ and is hence finite. For $\alpha\in\Xi_R$, we set the hyperplane
\[\alpha^\perp=\{\theta\in M_\mathcal{B}|\theta(\alpha)=0\}.\]
\begin{lemma}\label{polyhedral cell decomposition}
    Let $\Xi_R=\{\alpha_1,\dots,\alpha_m\}$. For
$s=(s_1,\dots,s_m)\in\{-1,0,1\}^m$, set
\[\Omega_s=\{\theta\in M_\mathcal{B}|\sgn\theta(\alpha_j)=s_j,j=1,\dots,m\}.\]
Then $\{\Omega_s|s\in\{-1,0,1\}^m,\Omega_s\neq\varnothing\}$ forms a finite rational polyhedral cell decomposition of $M_\mathcal{B}$. Moreover, each nonempty $\Omega_s$ is the relative interior of its closure $\overline{\Omega_s}$, which is a rational polyhedral cone. The three maps
\begin{align*}
    \mathcal{X}:M_\mathcal{B}&\rightarrow G_{\phi,R}\\
    \theta&\mapsto\1_{\mathcal{X}_\theta},
\end{align*}
where $\mathcal{X}\in\{\mathcal{C},\mathcal{T},\overline{\mathcal{T}}\}$, are constant on each $\Omega_s$.
\end{lemma}
We set
\[\mathcal{F}_R=\{\overline{\Omega_s}|s\in\{-1,0,1\}^m,\Omega_s\neq\varnothing\},\]
so each element in $\mathcal{F}_R$ is a cone. A cone $\mathfrak{d}\in\mathcal{F}_R$ is called a \emph{wall} if $\codim\mathfrak{d}=1$. It is easy to see that each wall is a subset of a hyperplane $\alpha^\perp$ where $\alpha\in\Xi_R$. We can pick a primitive vector $n_{\mathfrak{d}}\in\Gamma_\mathcal{B}$ parallel to $\alpha$ such that $h_\phi(n_{\mathfrak{d}})>0$. For any primitive vector $n$ with $h_\phi(n)>0$, we set
\[\mathfrak{g}_{n,R}=\bigoplus_{k\geqslant1,h_\phi(kn)\leqslant R}\mathcal{H}(\mathcal{B})_{kn},G_{n,R}=\exp_R\mathfrak{g}_{n,R}.\]

We write
\[\mathcal{A}_R=\bigcup_{\alpha\in\Xi_R}\alpha^\perp,\mathcal{W}_R=\{\mathfrak{d}\in\mathcal{F}_R|\codim\mathfrak{d}=1\},\mathcal{J}_R=\bigcup_{\codim\mathfrak{f}\geqslant2}\mathfrak{f}.\]

\begin{definition}
    A \emph{scattering diagram} is an assignment
    \[\mathfrak{D}=(g_{\mathfrak{d}})_{\mathfrak{d}\in\mathcal{W}_R},g_{\mathfrak{d}}\in G_{n_{\mathfrak{d}},R}.\]
    The \emph{support} of a scattering diagram $\mathfrak{D}$ is
    \[\supp\mathfrak{D}=\bigcup_{g_{\mathfrak{d}}\neq u_0}\mathfrak{d}.\]
\end{definition}
Now let $\gamma:[0,1]\rightarrow M_\mathcal{B}$ be a piecewise $C^1$ path. We say $\gamma$ is \emph{generic} if
\begin{itemize}
    \item $\gamma(0),\gamma(1)\notin\mathcal{A}_R$;
    \item $\gamma([0,1])\cap\mathcal{J}_R=\varnothing$;
    \item $\gamma$ meets $\mathcal{A}_R$ at only finitely many points $\gamma(t_i),i=1,\dots,k$ where $0<t_1<\dots<t_k<1$ and is $C^1$ near each $t_i$ with
    \[\left.\dfrac{\diff}{\diff t}\gamma(t)(n_{\mathfrak d_i})\right|_{t=t_i}\neq0,i=1,\dots,k.\]
\end{itemize}
At that time, we define the ordered product along the path $\gamma$ to be
\begin{equation}\label{scattering path product}
 \mathfrak p_{\gamma,R}=g_{\mathfrak d_1}^{\nu_1}*\cdots*g_{\mathfrak d_k}^{\nu_k},
\end{equation}
where
\[\nu_i=\sgn\left(\left.\dfrac{\diff}{\diff t}\gamma(t)(n_{\mathfrak d_i})\right|_{t=t_i}\right)\in\{\pm1\},i=1,\dots,k.\]
It is worth warning that our convention on the ordered product is different from that of Bridgeland \cite{Bri2017}. In fact, we have
\[\Phi_{\mathfrak D}(\gamma)=g_{\mathfrak d_k}^{-\nu_k}*\cdots*g_{\mathfrak d_1}^{-\nu_1}=\mathfrak p_{\gamma,R}^{-1}.\]
\begin{definition}
    A scattering diagram $\mathfrak{D}$ is said to be \emph{consistent} if, for any pair of generic paths $\gamma_1,\gamma_2$ in $M_\mathcal{B}$ with $\gamma_1(0)=\gamma_2(0)$ and $\gamma_1(1)=\gamma_2(1)$, we have $\mathfrak{p}_{\gamma_1,R}=\mathfrak{p}_{\gamma_2,R}$.
\end{definition}
We fix two special areas in $M_\mathcal{B}$:
\begin{align*}
    M_+&=\{\theta\in M_\mathcal{B}|\theta(\beta)>0\text{ for any }\beta\in\Lambda_\phi^+\},\\
    M_-&=\{\theta\in M_\mathcal{B}|\theta(\beta)<0\text{ for any }\beta\in\Lambda_\phi^+\}.
\end{align*}
Let $\mathfrak{D}$ be a consistent scattering diagram. We define 
\[\mathfrak{p}_{\mathfrak{D},R}=\mathfrak{p}_{\gamma,R},\]
where $\gamma$ is a generic path in $M_\mathcal{B}$ with $\gamma(0)\in M_+$ and $\gamma(1)\in M_-$. It is well-defined by the consistency of $\mathfrak{D}$.

We now construct a scattering diagram $\mathfrak{D}_{\sigma,\phi}^R$. For any $\mathfrak{d}\in\mathcal{W}_R$, we pick a $\theta$ in its relative interior and set $g_{\mathfrak{d}}=\1_{\mathcal{C}_\theta}\in G_{n_{\mathfrak{d}},R}$, where we always omit the projection map from $G_{\phi}$ to its quotients when there is no ambiguity. It is independent of the choice of $\theta$ by Lemma \ref{polyhedral cell decomposition}.
\begin{theorem}
    Let $\mathfrak{D}_{\sigma,\phi}^R$ be defined as above.
    \begin{enumerate}
        \item Let $\gamma:[0,1]\rightarrow M_\mathcal{B}$ be a generic path in $M_\mathcal{B}$. Then
        \[\mathfrak{p}_{\gamma,R}=\1_{\mathcal{T}_{\gamma(0)}}^{-1}*\1_{\mathcal{T}_{\gamma(1)}}.\]
        As a consequence, $\mathfrak{D}_{\sigma,\phi}^R$ is consistent.
        \item $\mathfrak{p}_{\mathfrak{D}_{\sigma,\phi}^R,R}=\1_R^{-1}$.
    \end{enumerate}
\end{theorem}
\begin{proof}
    Note that (2) follows immediately from (1) since $\mathcal{T}_\theta$ is $\mathcal{B}$ when $\theta\in M_+$ and is $0$ when $\theta\in M_-$. So we only need to prove (1). It suffices to prove the case that $\gamma$ meets only one wall $\mathfrak{d}\subseteq n^\perp$, where $n$ is primitive with $h_\phi(n)>0$, and the general case follows by induction.

    Fix $\eta\in M_\mathcal{B}$ such that $\eta(n)>0$. For any $\theta_0$ in the relative interior of $\mathfrak{d}$, set $\theta_\pm=\theta_0\pm\varepsilon\eta$, where $\varepsilon>0$ is small enough such that
    \[\varepsilon|\eta(\alpha)|<|\theta_0(\alpha)|,\alpha\in\Xi_R,\theta_0(\alpha)\neq0.\]
    What we need to show is
    \[g_{\mathfrak{d}}=\1_{\mathcal{C}_{\theta_0}}=\1_{\mathcal{T}_{\theta_-}}^{-1}*\1_{\mathcal{T}_{\theta_+}}.\]
    For any $\alpha\in\Xi_R$, we compare the signs of $\theta_*(\alpha),*\in\{0,-,+\}$.

    When $\theta_0(\alpha)=0$, there must exist a $k\in\mathbb{Z}_+$ such that $\alpha=kn$. Thus, we have
    \[\sgn\theta_*(\alpha)=\sgn(*\varepsilon\eta(\alpha))=*\sgn\eta(n)=*,*\in\{+,-\}.\]
    When $\theta_0(\alpha)\neq0$, the choice of $\varepsilon$ ensures that the signs of $\theta_*(\alpha)$ are the same. Therefore, we have the following table of signs:
    \[\begin{array}{c|c|c}
  \sgn\theta_0(\alpha)&\sgn\theta_-(\alpha)&\sgn\theta_+(\alpha)\\ \hline
  +&+&+\\
  0&-&+\\
  -&-&-
 \end{array}.\]
 Thus for any $E\in\mathcal{B}$ such that $m_\sigma(E)=h_\phi(\cl E)\leqslant R$, we have
 \[E\in\mathcal{T}_{\theta_-}\iff E\in\mathcal{T}_{\theta_0}, E\in\mathcal{T}_{\theta_+}\iff E\in\overline{\mathcal{T}}_{\theta_0},\]
 and hence
 \[\1_{\mathcal{T}_{\theta_-}}^{-1}*\1_{\mathcal{T}_{\theta_+}}=\1_{\mathcal{T}_{\theta_0}}^{-1}*\1_{\overline{\mathcal{T}}_{\theta_0}}=\1_{\mathcal{C}_{\theta_0}}\in G_{\phi,R},\]
 where the last equality follows from Corollary \ref{torsion semistable decomposition, hall}.
\end{proof}
\subsection{Nearby stability conditions as linear functionals}
For a \emph{center chart}, we mean a triple
\[c=(\sigma,\phi,R),\]
where $\sigma\in\Stab_\cl(\mathcal{D})$ satisfies the support property and the finite counting condition, $\phi\in\mathbb{R}$ and $R>0$. We write $V_c=V_R(\sigma)\subseteq\Stab_\cl(\mathcal{D})$.

Let $\tau=(Z_\tau,\mathcal{P}_\tau)\in V_c$ be a stability condition near $\sigma$. Given $\psi\in(\phi-\eta_R,\phi+\eta_R)$, we set
\begin{align*}
    \Theta(\tau,\psi):\Gamma_\mathcal{B}&\rightarrow\mathbb{R}\\
    \beta&\mapsto\im\left(e^{-i\pi\psi}Z_\tau(\beta)\right).
\end{align*}
One can check that $\Theta(\tau,\psi)\in M_\mathcal{B}$. Moreover, if $0\neq E\in\mathcal{B}$ with $m_\sigma(E)=h_\phi(\cl E)\leqslant R$ satisfies that $Z_\tau(\cl E)=r_Ee^{i\pi\alpha_E}$ for some $r_E>0$ and $\alpha_E\in(\phi-\eta_R,\phi+\eta_R)$, then we have
\[\Theta(\tau,\psi)(\cl E)=r_E\sin(\pi(\alpha_E-\psi)).\]
\begin{lemma}\label{Ctheta=Ptau}
    Let $\theta=\Theta(\tau,\psi)$. Then for any $0\neq E\in\mathcal{B}$ with $m_\sigma(E)\leqslant R$, we have $E\in\mathcal{C}_\theta$ if and only if $E\in\mathcal{P}_\tau(\psi)$.
\end{lemma}
\begin{proof}
    By the definition of $\theta=\Theta(\tau,\psi)$, $E\in\mathcal{C}_\theta$ if and only if $\alpha_E=\psi$ and $\alpha_A\leqslant\psi$ for any $0\neq A\subseteq E$. Note that these coincide with the condition for $E$ to be $Z_\tau|_\mathcal{B}$-semistable with phase $\psi$, so the result follows from Lemma \ref{semistable agree}.
\end{proof}
The following proposition is the translation of Lemma \ref{Ctheta=Ptau} into the language of Hall algebras, by applying Proposition \ref{nearby} and Corollary \ref{delta agree}.
\begin{proposition}\label{delta and epsilon}
    Let $\beta\in\Lambda_{\phi,R}$, and let $\psi\in(\phi-\eta_R,\phi+\eta_R)$ be the local phase of $Z_\tau(\beta)$. Set $\theta=\Theta(\tau,\psi)$. Then
    \begin{enumerate}
        \item $(\1_{\mathcal{C}_\theta})_\beta=\delta^\tau_{\psi,\beta}$;
        \item If we define $\epsilon^{c,\tau}_{\psi,\beta}:= (\log\1_{\mathcal{C}_\theta})_\beta$, then
        \[\epsilon^{c,\tau}_{\psi,\beta}=\epsilon^\tau_{\psi,\beta}.\]
    \end{enumerate}
\end{proposition}
\begin{corollary}\label{switching}
    Let $c=(\sigma,\phi,R)$ and $d=(\sigma',\phi',R')$ be two center charts. Suppose that $\tau\in V_c\cap V_d,\beta\in\Lambda^c_{\phi,R}\cap\Lambda^d_{\phi',R'}$ with $Z_\tau(\beta)\in\mathbb{R}_+e^{i\pi\psi}$ where $\psi\in(\phi-\eta_R,\phi+\eta_R)\cap(\phi'-\eta'_{R'},\phi'+\eta'_{R'})$. Then
    \[\epsilon^{c,\tau}_{\psi,\beta}=\epsilon^{\tau}_{\psi,\beta}=\epsilon^{d,\tau}_{\psi,\beta}.\]
\end{corollary}
\begin{proof}
    We obtain the result immediately by applying Proposition \ref{delta and epsilon}(2) to the two center charts respectively.
\end{proof}
Let $c=(\sigma,\phi,R)$ be a center chart. For $\tau\in V_c$, we define the \emph{active homogeneous frame} of $c$ at $\tau$ to be the finite set
\[\mathscr{E}_c(\tau)=\left\{\epsilon^\tau_{\phi_\tau(\beta),\beta}\middle|0\neq\beta\in\Lambda_{\phi,R},\epsilon^\tau_{\phi_\tau(\beta),\beta}\neq0\right\}\subseteq\mathcal{H}(\mathcal{D}).\]
Let $P$ be a finite Lie polynomial in formal symbols $x_{\psi,\beta},\psi\in\mathbb{R},\beta\in\Gamma$. We set
\[\ev_\tau(x_{\psi,\beta})=\epsilon^\tau_{\psi,\beta},\]
and extend it to a Lie homomorphism
\[\ev_\tau:\Lie_{\mathbb{Q}}\{x_{\psi,\beta}|\psi\in\mathbb{R},\beta\in\Gamma\}\rightarrow\mathcal{H}(\mathcal{D}).\]

\begin{theorem}
    Let $c=(\sigma,\phi,R)$ and $d=(\sigma',\phi',R')$ be two center charts and $\tau_0,\tau_1\in V_c\cap V_d$. Let $P_c$ and $P_d$ be Lie polynomials in $\epsilon^{\tau_0}_{\phi_{\tau_0}(\beta),\beta}$'s to express $\epsilon^{\tau_1}_{\psi,\gamma}$ in the sense of Theorem \ref{Joyce transport} in $V_c$ and $V_d$, respectively. Then
    \[\ev_{\tau_0}(P_c)=\epsilon^{\tau_1}_{\psi,\gamma}=\ev_{\tau_0}(P_d),\]
    although the Lie polynomials $P_c$ and $P_d$ may be different.
\end{theorem}
\begin{proof}
    The result follows immediately from Theorem \ref{Joyce transport} and Corollary \ref{switching}.
\end{proof}
\subsection{Derived scattering groupoids}
Let $\Sigma$ be a connected component of $\Stab_\cl(\mathcal{D})$ such that there exists a $\sigma\in\Sigma$ satisfying the support property and the finite counting condition. Lemma \ref{one point to the component} ensures that every stability condition in $\Sigma$ satisfies these two conditions.
\begin{definition}
    A \emph{finite scattering frame} at $\tau\in\Sigma$ is a finite labeled family
    \[\mathscr{F}=\left(\tau;\{\epsilon^\tau_{\psi_i,\gamma_i}|i\in I\}\right),\]
    where $\epsilon^\tau_{\psi_i,\gamma_i}\neq0$ for each $i\in I$.
\end{definition}
Recall that a \emph{tree} is a finite connected graph without unoriented cycles. For a tree $T=(T_0,T
_1)$ with $|T_0|\geqslant2$, let $p$ be a vertex of $T$. Write $c_p$ for the number of edges in $T_1$ having $p$ as an endpoint. Let
\[L=\{p\in T_0|c_p=1\}.\]
Fix a vertex $r\in T_0$. We call it the \emph{root} of $T$. Then a vertex $r\neq l\in L$ is called a \emph{leaf}. A vertex $r\neq p\in T_0-L$ is called a \emph{non-leaf vertex}. In the special case that $T_0=\{p\}$, we call $p$ both the root and the leaf of $T$, and $T$ has no non-leaf vertex in this case.

For any vertex $r\neq p\in T_0$, there is a sequence of edges $e_1,\dots,e_m$ forming the unique simple path in $T$ from $p$ to $r$. At that time, $p$ is an endpoint of $e_1$. We call the other endpoint $q$ the \emph{parent vertex} of $p$ and call $p$ a \emph{child vertex} of $q$. We call $p$ and $q$ the \emph{child end} and the \emph{parent end} of $e_1$, respectively.
\begin{definition}
    Let $\tau_0,\tau_1\in\Sigma$ and $\mathscr{F}_j$ be a finite scattering frame at $\tau_j$ for each $j=0,1$. A \emph{transport tree} from $\mathscr{F}_0$ to $\mathscr{F}_1$ is a tree together with the following data:
    \begin{itemize}
        \item The root of the tree is labelled by a target element $\epsilon^{\tau_1}_{\psi,\gamma}$ in $\mathscr{F}_1$ together with a Lie polynomial.
        \item Each leaf of the tree is labelled by a source element $\epsilon^{\tau_0}_{\psi,\gamma}$ in $\mathscr{F}_0$.
        \item Each non-leaf vertex of the tree is labelled by an element $\epsilon^\tau_{\psi,\gamma}$ together with a Lie polynomial such that, for any vertex of the tree that is not a leaf, the element labelling it is the evaluation of the Lie polynomial at all elements labelling its child vertices. 
    \end{itemize}
\end{definition}
Figure \ref{fig:scattering-substitution-tree} illustrates this process schematically. The polynomials $P,Q_1,Q_2$ stand for local substitutions, and the resulting source
expression is
\[Y=P(Q_1(X_1,X_2),Q_2(X_3,X_4)).\]
The downward arrows record expansion; evaluation runs
upwards after the source values have been inserted.

\par\medskip
\noindent\begin{minipage}{\linewidth}
\centering
\begin{tikzpicture}[
  >=Latex,
  every node/.style={font=\small},
  formula/.style={draw,rounded corners=2pt,fill=black!3,
                  inner xsep=8pt,inner ysep=7pt},
  source/.style={draw,rounded corners=2pt,
                 minimum width=14mm,minimum height=9mm},
  edge/.style={->,semithick}
]
 \node[formula] (Y) at (0,0) {$Y=P(A_1,A_2)$};
 \node[formula] (A1) at (-2.8,-1.45)
     {$A_1=Q_1(X_1,X_2)$};
 \node[formula] (A2) at (2.8,-1.45)
     {$A_2=Q_2(X_3,X_4)$};
 \node[source] (X1) at (-4.15,-2.95) {$X_1$};
 \node[source] (X2) at (-1.45,-2.95) {$X_2$};
 \node[source] (X3) at (1.45,-2.95) {$X_3$};
 \node[source] (X4) at (4.15,-2.95) {$X_4$};
 \draw[edge] (Y) -- (A1);
 \draw[edge] (Y) -- (A2);
 \draw[edge] (A1) -- (X1);
 \draw[edge] (A1) -- (X2);
 \draw[edge] (A2) -- (X3);
 \draw[edge] (A2) -- (X4);
 \node[above=3pt] at (Y.north) {target coordinate};
 \node at (0,-3.65) {source coordinates};
 \draw[edge] (-5.45,-0.35) -- (-5.45,-2.7);
 \node[rotate=90] at (-5.85,-1.52) {expansion};
 \draw[edge] (5.45,-2.7) -- (5.45,-0.35);
 \node[rotate=90] at (5.85,-1.52) {evaluation};
\end{tikzpicture}

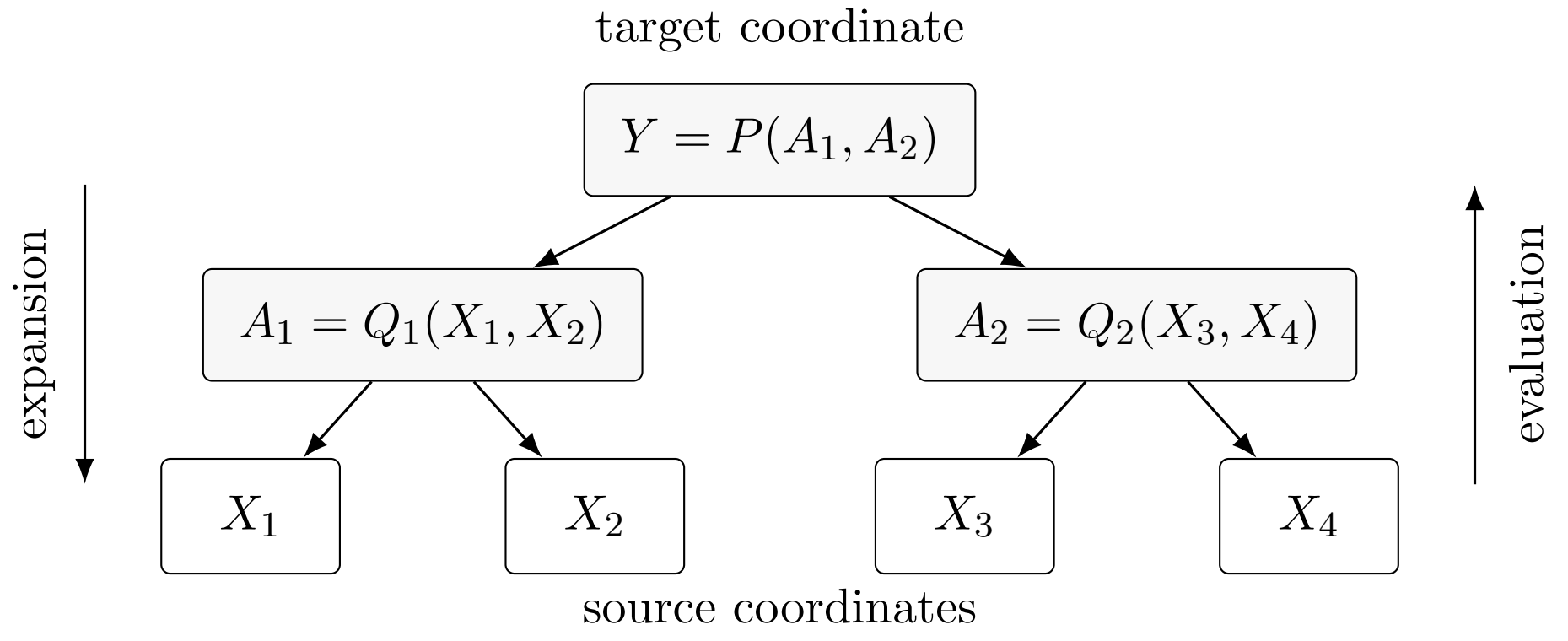
\captionof{figure}{A finite substitution tree.}
\label{fig:scattering-substitution-tree}
\end{minipage}
\par\medskip

By Theorem \ref{Joyce transport}, when the source frame $\mathscr{F}_0$ has enough source elements, we can construct a transport tree from $\mathscr{F}_0$ to $\mathscr{F}_1$. We call a transport tree from $\mathscr{F}_0$ to $\mathscr{F}_1$ \emph{admissible} if the following conditions hold:
\begin{itemize}
    \item For each vertex that is not a leaf, the expression of its element as a Lie polynomial in the elements on its child vertices is obtained in a center chart $c=(\sigma,\phi,R)$ by Theorem \ref{local Joyce transformation}.
    \item For any non-leaf vertex $p$, let $\epsilon^\tau_{\psi,\gamma}$ be the element on it. Let $c,d$ be two center charts such that the expression of $\epsilon^\tau_{\psi,\gamma}$ in the elements on the child vertices of $p$ is in $c$ and the expression of the element on the parent vertex $q$ of $p$ in the elements on the child vertices of $q$ is in $d$. Then $\tau\in V_c\cap V_d$. At that time, the element $\epsilon^{c,\tau}_{\psi,\gamma}=\epsilon^\tau_{\psi,\gamma}=\epsilon^{d,\tau}_{\psi,\gamma}$ is identified by Corollary \ref{switching}. 
\end{itemize}
\begin{definition}
    An \emph{admissible transport forest} $\mathscr{T}$ from $\mathscr{F}_0$ to $\mathscr{F}_1$ is a finite set of admissible transport trees from $\mathscr{F}_0$ to $\mathscr{F}_1$ such that every target element $\epsilon^{\tau_1}_{\psi_i,\gamma_i}$ occurs and only occurs in the root of one transport tree in $\mathscr{T}$.
\end{definition}
    
Let $\mathscr{F}_0=(\tau_0;\{X_i\}_{i\in I}),\mathscr{F}_1=(\tau_1;\{Y_j\}_{j\in J})$ and $\mathscr{T}=\{T_j|j\in J\},\mathscr{T}'=\{T'_j|j\in J\}$ be two admissible transport forests from $\mathscr{F}_0$ to $\mathscr{F}_1$. Let $P_j$ and $P'_j$ be the resulting Lie polynomials in $T_j$ and $T'_j$ to express $Y_j$  in $X_i,i\in I$, respectively. We say that $\mathscr{T}$ and $\mathscr{T}'$ are \emph{equivalent} if
\[\ev_{\tau_0}(P_j)=\ev_{\tau_0}(P'_j),j\in J.\]
In fact, this equivalence relation is trivial, since all admissible transport forests from $\mathscr{F}_0$ to $\mathscr{F}_1$ are equivalent.
\begin{remark}
    We retain this quotient notation because the same construction can be used with a finer equivalence relation. With the present equivalence relation, the transport category is thin.
\end{remark}
\begin{definition}
    We define the \emph{transport category} $\mathfrak{Tr}(\mathcal{D},\Sigma)$ on $\Sigma$ with the following data:
    \begin{itemize}
        \item The objects are finite scattering frames on $\tau$, where $\tau\in\Sigma$.
        \item A morphism from $\mathscr{F}_0$ to $\mathscr{F}_1$ is an equivalence class of admissible transport forests from $\mathscr{F}_0$ to $\mathscr{F}_1$.
        \item For each $\mathscr{F}\in\mathfrak{Tr}(\mathcal{D},\Sigma)$, the identity morphism at $\mathscr{F}$ is defined to be the equivalence class of the forest consisting of one-point trees for each Joyce element in $\mathscr{F}$.
        \item The composition of two forests $\mathscr{T}:\mathscr{F}_0=(\sigma_0;\{X_i\}_{i\in I})\rightarrow\mathscr{F}_1=(\sigma_1;\{Y_j\}_{j\in J})$ and $\mathscr{T}':\mathscr{F}_1\rightarrow\mathscr{F}_2=(\sigma_2;\{Z_k\}_{k\in K})$ is defined to be the forest $\mathscr{T}'\circ\mathscr{T}=\{(T'\circ T)_k|k\in K\}$ such that the tree $(T'\circ T)_k$ with each Joyce element $Z_k$ in $\mathscr{F}_2$ is obtained by the following process: Pick the tree $T'_k$ in $\mathscr{T}'$ with $Z_k$ on its root and find the elements on its leaves $Y_{j_1},\dots,Y_{j_s}$ in $\mathscr{F}_1$, which are in the frame $\mathscr{F}_1$. Then pick the trees $T_{j_1},\dots,T_{j_s}\in\mathscr{T}$ such that $T_{j_t}$ is the tree in $\mathscr{T}$ with $Y_{j_t}$ on its root, $t=1,\dots,s$. Now let $(T'\circ T)_k$ be the tree obtained by gluing the roots of $T_{j_1},\dots,T_{j_s}$ to the leaves labelled by $Y_{j_1},\dots,Y_{j_s}$ of $T'_k$ respectively. One can easily check that $(T'\circ T)_k$ is indeed a tree and is admissible using Corollary \ref{switching}. The composition in $\mathfrak{Tr}(\mathcal{D},\Sigma)$ is induced by the composition of admissible transport forests.
    \end{itemize}
\end{definition}
Since the equivalence relation is trivial, for any $\mathscr{F}_0,\mathscr{F}_1\in\mathfrak{Tr}(\mathcal{D},\Sigma)$, we have
\[\Hom_{\mathfrak{Tr}(\mathcal{D},\Sigma)}(\mathscr{F}_0,\mathscr{F}_1)=\{*\}\text{ or }\varnothing,\]
depending on whether an admissible transport forest exists. So one can check that the composition defined above is well-defined and satisfies the associative law, and hence $\mathfrak{Tr}(\mathcal{D},\Sigma)$ is indeed a category. If there exist two admissible transport forests $\mathscr{T}:\mathscr{F}_0\rightarrow\mathscr{F}_1$ and $\mathscr{T}':\mathscr{F}_1\rightarrow\mathscr{F}_0$,
then their isomorphism classes must be the inverses of each other. Now we define
\begin{definition}
    Let $\mathfrak{Scat}(\mathcal{D},\Sigma)$ be the following category: It has the same objects as $\mathfrak{Tr}(\mathcal{D},\Sigma)$, and its morphisms are precisely the isomorphisms in $\mathfrak{Tr}(\mathcal{D},\Sigma)$. $\mathfrak{Scat}(\mathcal{D},\Sigma)$ is a groupoid and is called the \emph{scattering transport groupoid} of $\Sigma$.
\end{definition}
\begin{theorem}
    Let $\sigma_0,\sigma_1\in\Sigma$ and $\sigma:[0,1]\rightarrow\Sigma$ be a piecewise $C^1$ path from $\sigma_0$ to $\sigma_1$. Choose $R>0$ as in the proof of Theorem \ref{Joyce transport}. Then the proof process of Theorem \ref{Joyce transport} gives a construction of an admissible transport forest from a finite scattering frame at $\sigma_0$ to one at $\sigma_1$. In particular, for any finite scattering frame $\mathscr{F}_1$ at $\sigma_1$, there is a finite scattering frame $\mathscr{F}_0$ at $\sigma_0$, such that there exists an admissible transport forest from $\mathscr{F}_0$ to $\mathscr{F}_1$.
\end{theorem}
\section{The finite-length heart case}
In this section, we assume that $\mathcal{D}$ has a heart of finite length with only finitely many simple objects up to isomorphism. For example, $\mathcal{D}=D^b(Q)$, where $Q$ is a finite acyclic quiver. We will give a finite length version of Theorem \ref{Joyce transport} and provide an alternative proof in the special case.

Let $\mathcal{A}$ be such a heart of $\mathcal{D}$. Then we have
\[K_0(\mathcal{D})\cong K_0(\mathcal{A})=\mathbb{Z}^n,n\in\mathbb{Z}_+.\]
So we will fix $\Gamma=K_0(\mathcal{D})$ and $\cl=\id_{K_0(\mathcal{D})}$, and therefore $\Stab_\cl(\mathcal{D})=\Stab(\mathcal{D})$.
\subsection{Rotated algebraic stability conditions}
The following notion follows Qiu and Woolf \cite{QW2018}.
\begin{definition}
    We say that a stability condition $\sigma=(Z,\mathcal{P})$ on $\mathcal{D}$ is \emph{algebraic} if its heart $\mathcal{P}(0,1]$ is of finite length with only finitely many isomorphism classes of simple objects. We call it \emph{rotated algebraic} if it has a rotated heart $\mathcal{P}(a,a+1]$ that is of finite length with only finitely many isomorphism classes of simple objects.
\end{definition}
We write $\Stab_\fl(\mathcal{D})$ for the subset of $\Stab(\mathcal{D})$ consisting of rotated algebraic stability conditions, where ``fl'' stands for ``finite length''.
\begin{lemma}[{\cite[Lemma 3.1; Lemma 3.11]{QW2018}}]\label{fl neighbourhood}
    Let $\sigma\in\Stab_\fl(\mathcal{D})$ with a rotated heart $\mathcal{A}$ of finite length with only finitely many isomorphism classes of simple objects. Then there exist  $a\in\mathbb{R}$ and $\eta>0$, such that for any $\tau\in V_\eta(\sigma)=\{\tau\in\Stab(\mathcal{D})|d(\sigma,\tau)<\eta\}$, we have
    \[\mathcal{P}_\tau(a,a+1]=\mathcal{A}.\]
\end{lemma}
\begin{proof}
    Let $S_1,\dots,S_n$ be the simple objects in $\mathcal{A}$, and set
    \[b=\min_{1\leqslant j\leqslant n}\phi_\sigma^-(S_j),c=\max_{1\leqslant j\leqslant n}\phi_\sigma^+(S_j),\]
    so we must have $0\leqslant c-b<1$. Thus, we can choose $a\in\mathbb{R},\eta>0$, such that
    \[a+2\eta<b\leqslant c<a+1-2\eta.\]
    So $\mathcal{A}=\mathcal{P}(a,a+1]$. When $d(\sigma,\tau)<\eta$, since $|\phi_\sigma^\pm(S_j)-\phi_\tau^\pm(S_j)|<\eta,j=1,\dots,n$, we have $S_1,\dots,S_n\in\mathcal{P}_\tau(a,a+1]$ and hence $\mathcal{A}\subseteq\mathcal{P}_\tau(a,a+1]$. Conversely, for any $0\neq E\in\mathcal{P}_\tau(a,a+1]$, we have $\phi_\sigma^\pm(E)\in(a-\eta,a+1+\eta)$. By the assumption of $\sigma$, for any $0\neq E\in\mathcal{D}$, $\phi_\sigma^\pm(E)$ can only lie in $[b+k,c+k],k\in\mathbb{Z}$. By the choice of $a$ and $\eta$, we must have $\phi_\sigma^\pm(E)\in[b,c]\subseteq(a,a+1]$, so $E\in\mathcal{P}_\sigma(a,a+1]=\mathcal{A}$, which proves that $\mathcal{P}_\tau(a,a+1]\subseteq\mathcal{A}.$
\end{proof}
\subsection{Calculate in a finite-length heart}
Let $\mathcal{A}$ be a finite length heart of $\mathcal{D}$ with only finitely many isomorphism classes of simple objects. By \cite[Proposition 6.1]{Toe2006}, we have the following lemma:
\begin{lemma}
    The subalgebra $\mathcal{H}(\mathcal{A})\subseteq\mathcal{H}(\mathcal{D})$ generated by $u_A,A\in\mathcal{A}$ coincides with the Ringel-Hall algebra of $\mathcal{A}$.
\end{lemma}

Let $\mathcal{A}$ be a finite length heart of $\mathcal{D}$ with only finitely many isomorphism classes of simple objects and let $Z$ be a rotated stability function on it. For any $0\neq\gamma\in K_0(\mathcal{A})=K_0(\mathcal{D})=\mathbb{Z}^n$, we can define
\[\kappa^\mathcal{A}_\gamma=\1_{\{X\in\Iso(\mathcal{A})|\udv X=\gamma\}},\delta^Z_\gamma=\1_{\{X\in\Iso(\mathcal{A})|\udv X=\gamma,X\text{ is $Z$-semistable}\}}\in\mathcal{H}(\mathcal{A})\]
as in \cite{BT2012}. The corresponding Joyce element is
\[\epsilon^Z_\gamma=\sum_{\begin{subarray}{c}
    \\ \gamma_1+\dots+\gamma_r=\gamma\\
   \phi_Z(\gamma_1)=\dots=\phi_Z(\gamma_r)=\phi_Z(\gamma) 
\end{subarray}}\frac{(-1)^{r-1}}{r}\delta^Z_{\gamma_1}*\dots*\delta^Z_{\gamma_r}.\]

For $\sigma=(Z,\mathcal{P})\in\Stab_\fl(\mathcal{D})$ and a finite length rotated heart $\mathcal{A}=\mathcal{P}(a,a+1]$ with only finitely many isomorphism classes of simple objects, we know $Z|_\mathcal{A}$ is a rotated stability function on $\mathcal{A}$ with HN property. Using the notation in Section 4, one can easily check that
\[\kappa^\mathcal{A}_\gamma=\kappa^\sigma_{(a,a+1],\gamma},\delta^{Z|_\mathcal{A}}_\gamma=\delta^\sigma_{\phi(\gamma),\gamma}, \epsilon^{Z|_\mathcal{A}}_{\gamma}=\epsilon^\sigma_{\phi(\gamma),\gamma}\]
where $\phi(\gamma)\in(a,a+1]$ is the phase of $Z(\gamma)$ in one branch, and the formulas (\ref{HN identity in sector})-(\ref{Joyce inversion}) still hold. Note that $(a,a+1]$ is not a strict sector, but we can also define $\kappa^\sigma_{(a,a+1],\gamma}$ with the same form, since $\mathcal{P}_\sigma(a,a+1]=\mathcal{P}_\sigma[b,c]$ for some $[b,c]\subsetneqq(a,a+1]$ by the proof of Lemma \ref{fl neighbourhood} and $[b,c]$ must be a strict sector. By the same argument as in the proof of Proposition \ref{decomposition of 1}, we obtain

\begin{proposition}\label{decomposition of 1,fl}
    \[\1_\mathcal{A}=\prod_{\phi\in(a,a+1]}^\leftarrow\exp\left(\sum_{\phi(\gamma)=\phi}\epsilon^\sigma_{\phi,\gamma}\right)\in\hat{\mathcal{H}}(\mathcal{A}).\]
\end{proposition}
\subsection{Finite length Joyce transport}
Let $\sigma\in\Stab_\fl(\mathcal{D})$. We set $\mathfrak{L}_\sigma$ to be the Lie subalgebra of $\mathcal{H}(\mathcal{D})$ generated by the subset
\[\{\epsilon^\sigma_{\phi,\gamma}|\phi\in\mathbb{R},\gamma\in K_0(\mathcal{D})=\mathbb{Z}^n\}.\]
The following theorem is the finite-length version of Theorem \ref{Joyce transport}.
\begin{theorem}[Finite length Joyce transport]\label{finite length Joyce transport}
    Let $\sigma:[0,1]\rightarrow\Stab_\fl(\mathcal{D})$ be a continuous path. For any $\phi\in\mathbb{R},\gamma\in K_0(\mathcal{D})=\mathbb{Z}^n$, the Joyce element $\epsilon^{\sigma_1}_{\phi,\gamma}$ can be written as a finite Lie polynomial in $\epsilon^{\sigma_0}_{\psi,\beta},\psi\in\mathbb{R},\beta\in K_0(\mathcal{D})=\mathbb{Z}^n$ with coefficients in $\mathbb{Q}$, with each Lie monomial having total degree $\gamma$. As a consequence, we have
    \[\mathfrak{L}_{\sigma_0}=\mathfrak{L}_{\sigma_1}.\]
\end{theorem}
\begin{proof}
    For any $t\in[0,1]$, we can pick an open segment $J_t\subseteq\mathbb{R}$ such that, for any $s\in J_t\cap[0,1]$, we have $\sigma_s\in V_{\eta_t}(\sigma_t)$ for some $\eta_t>0$ in Lemma \ref{fl neighbourhood}.
    By the compactness of $[0,1]$, one can choose a finite subset $\{t_1,\dots,t_N\}\subseteq[0,1]$ such that
    \[[0,1]\subseteq\bigcup_{j=1}^N J_{t_j}=\bigcup_{j=1}^N(c_{t_j},d_{t_j}).\]
    We can shrink each segment $(c_{t_j},d_{t_j})$ such that $t_j<c_{t_{j+1}}<d_{t_j}<t_{j+1},j=1,\dots,N-1$. Then we can pick $0=b_0<b_1<\dots<b_N=1$ such that $t_j\in[b_{j-1},b_j]\subseteq J_{t_j}$ for each $j=1,\dots,N$. Thus, for any $t\in[b_{j-1},b_j],j=1,\dots,N$, there exists an $a_j\in\mathbb{R}$ and a finite length heart $\mathcal{A}_j$ of $\mathcal{D}$ with only finitely many isomorphism classes of simple objects, such that $\mathcal{P}_t(a_j,a_j+1]=\mathcal{A}_j$.

    Now consider each segment $[b_{j-1},b_j]$. It suffices to consider the Joyce element $\epsilon^{\sigma_{b_j}}_{\psi,\gamma}$ where $\psi\in(a_j,a_j+1]$, since if $\psi\notin(a_j,a_j+1]$ we can apply a shift. The stability conditions $\sigma_{b_{j-1}}$ and $\sigma_{b_j}$ have the common rotated heart $\mathcal{A}_j$. Let us work in it. By proposition \ref{decomposition of 1,fl}, we have
    \[\prod_{\phi\in(a_j,a_j+1]}^\leftarrow\exp\left(\sum_{\phi_{\sigma_{b_{j-1}}}(\gamma)=\phi}\epsilon^{\sigma_{b_{j-1}}}_{\phi,\gamma}\right)=\1_{\mathcal{A}_j}=\prod_{\phi\in(a_j,a_j+1]}^\leftarrow\exp\left(\sum_{\phi_{\sigma_{b_j}}(\gamma)=\phi}\epsilon^{\sigma_{b_j}}_{\phi,\gamma}\right)\in\hat{\mathcal{H}}(\mathcal{A}_j).\]
    By Theorem \ref{universal Lie}, $\epsilon^{\sigma_{b_j}}_{\phi,\gamma}$ can be expressed as a Lie polynomial in $\epsilon^{\sigma_{b_{j-1}}}_{\psi,\beta}$'s. The result follows by evaluation backwards along the path.
\end{proof}
\subsection{An application for Dynkin quivers}
Let $Q=(Q_0,Q_1,s,t)$ be a Dynkin quiver of ADE type with $n$ vertices. We write $Q_0=\{1,\dots,n\}$ so that every arrow $(\alpha:s(\alpha)\rightarrow t(\alpha))\in Q_1$ satisfies $s(\alpha)<t(\alpha)$. Let $\mathcal{D}_Q=D^b(\rep_kQ)$. Let $\Phi_+(Q)$ be the set of positive roots corresponding to the underlying graph of $Q$.

We quote the following two lemmas and omit their proofs.
\begin{lemma}[{\cite[Lemma 4.11]{QW2018}}]\label{Dynkin length}
Every heart of a stability condition on $\mathcal{D}_Q$ has finite length, with exactly $n$ isomorphism classes of simple objects, and at most $|\Phi_+(Q)|$ isomorphism classes of indecomposable objects. In particular,
\[\Stab_{\fl}(\mathcal{D}_Q)=\Stab(\mathcal{D}_Q).\]
\end{lemma}
\begin{lemma}[{\cite[Theorem 4.7]{Qiu2015}}]
\label{Dynkin path-connectedness}
$\Stab(\mathcal{D}_Q)$ is simply connected, and hence path connected.
\end{lemma}

Let $\sigma_*=(Z_*,\mathcal{P}_*)\in\Stab(\mathcal{D}_Q)$ be given by the heart $\mathcal{A}_Q=\rep_{\mathbb{F}_q}Q$ with
\[0<\phi_*(S_1)<\cdots<\phi_*(S_n)<1,\]
in the sense of Lemma \ref{heart}. The following lemma follows by calculating in the heart $\mathcal{A}_Q$.

\begin{lemma}
\label{standard stability condition}
We have
\[\mathcal{P}_*(\phi_*(S_i)+b)=\add(S_i[b])
 ,i\in Q_0,b\in\mathbb{Z},\]
and $\mathcal{P}_*(\phi)=0$ for all other $\phi\in\mathbb{R}$. Equivalently, all $\sigma_*$-semistable objects are precisely $S_i^{\oplus r}[b]$ with $i\in Q_0$, $r\in\mathbb{Z}_+$, and $b\in\mathbb{Z}$.
\end{lemma}
\begin{proof}
    It suffices to work in the heart $\mathcal{A}_Q$. Let $0\neq X\in\mathcal{A}_Q$. For $s=1,\dots,n+1$, set $F_sX$ to be the subrepresentation of $X$ with spaces $X_t$ at the vertices $t=s,s+1,\dots,n$ and $0$ on other vertices. The linear maps on the arrows are the corresponding restrictions. We have a filtration of $X$:
    \[0=F_{n+1}X\subseteq F_nX\subseteq\dots\subseteq F_1X=X,F_sX/F_{s+1}X=S_s^{\oplus\dim X_s}.\]
    Let $1\leqslant j_1<\dots<j_t\leqslant n$ such that $\{j_1,\dots,j_t\}=\{1\leqslant s\leqslant n|\dim X_s\neq0\}$. Then
    \[0\subsetneqq F_{j_t}X\subsetneqq F_{j_{t-1}}X\subsetneqq\dots\subsetneqq F_{j_1}X=X\]
    is the HN filtration of $X$ with respect to $Z_*|_{\mathcal{A}_Q}$. Therefore, $X$ is $Z_*|_{\mathcal{A}_Q}$-semistable if and only if $t=1$, and in this case $X=S_{j_1}^{\oplus\dim X_{j_1}}$.
\end{proof}
We call $\sigma_*$ a \emph{standard stability condition} on $\mathcal{D}_Q$. One can easily calculate that
\[\delta^{\sigma_*}_{\phi_*(S_i)+b,(-1)^b re_i}=u_{S_i^{\oplus r}[b]}\]
are the only nonzero $\delta$ functions, and hence the corresponding Joyce elements are
\begin{align*}
    \epsilon^{\sigma_*}_{\phi_*(S_i)+b,(-1)^b re_i}&=\sum_{r_1+\dots+r_k=r}\dfrac{(-1)^{k-1}}{k}\delta^{\sigma_*}_{\phi_*(S_i)+b,(-1)^b r_1e_i}*\dots*\delta^{\sigma_*}_{\phi_*(S_i)+b,(-1)^b r_ke_i}\\
    &=\sum_{r_1+\dots+r_k=r}\dfrac{(-1)^{k-1}}{k}u_{S_i^{\oplus r_1}[b]}*\dots*u_{S_i^{\oplus r_k}[b]}.
\end{align*}
We have
\[u_{S_i^{\oplus r}}*u_{S_i^{\oplus s}}=\left[\begin{array}{c}
     r+s  \\
     s 
\end{array}\right]_q u_{S_i^{\oplus(r+s)}},\]
where
\[\left[\begin{array}{c}
     r+s  \\
     s 
\end{array}\right]_q=\frac{\displaystyle{\prod_{j=1}^{r+s}(q^j-1)}}{\displaystyle{\prod_{j=1}^{r}(q^j-1)\prod_{j=1}^{s}(q^j-1)}}\]
is the Gaussian binomial coefficient. So
\begin{align*}
    \epsilon^{\sigma_*}_{\phi_*(S_i)+b,(-1)^b re_i}
    &=\sum_{r_1+\dots+r_k=r}\dfrac{(-1)^{k-1}}{k}u_{S_i^{\oplus r_1}[b]}*\dots*u_{S_i^{\oplus r_k}[b]}\\
    &=\sum_{r_1+\dots+r_k=r}\dfrac{(-1)^{k-1}}{k}\left[\begin{array}{c}
     r_1+r_2  \\
     r_1 
\end{array}\right]_q\dots\left[\begin{array}{c}
     r_1+\dots+r_k  \\
     r_1+\dots+r_{k-1}
\end{array}\right]_q u_{S_i^{\oplus r}[b]}\\
&=\sum_{r_1+\dots+r_k=r}\dfrac{(-1)^{k-1}}{k}\frac{\displaystyle{\prod_{j=1}^{r}(q^j-1)}}{\displaystyle{\prod_{j=1}^{r_1}(q^j-1)\dots\prod_{j=1}^{r_k}(q^j-1)}}u_{S_i^{\oplus r}[b]}\\
&=\dfrac{1}{r}\prod_{j=1}^{r-1}(1-q^j)u_{S_i^{\oplus r}[b]},
\end{align*}
where the last equality is based on the following combinatorial lemma.
\begin{lemma}
    Let $q,r\in\mathbb{Z}_+$. We have
    \[\sum_{r_1+\dots+r_k=r}\dfrac{(-1)^{k-1}}{k}\frac{\displaystyle{\prod_{j=1}^{r}(q^j-1)}}{\displaystyle{\prod_{j=1}^{r_1}(q^j-1)\dots\prod_{j=1}^{r_k}(q^j-1)}}=\dfrac{1}{r}\prod_{j=1}^{r-1}(1-q^j).\]
\end{lemma}
\begin{proof}
    Set
    \[d_k=\prod_{j=1}^k(1-q^j),k>0\]
    and $d_0=1$. Then we have $d_k=(1-q^k)d_{k-1},k>0$. Consider the power series
    \[F(x)=\sum_{k=0}^\infty\dfrac{x^k}{d_k}.\]
    We have
    \begin{align*}
        F(qx)&=\sum_{k=0}^\infty\dfrac{q^kx^k}{d_k}\\
        &=\sum_{k=0}^\infty\dfrac{x^k}{d_k}-\sum_{k=0}^\infty\dfrac{(1-q^k)x^k}{d_k}\\
        &=F(x)-\sum_{k=1}^\infty\dfrac{x^k}{d_{k-1}}\\
        &=F(x)-xF(x)\\
        &=(1-x)F(x).
    \end{align*}
    Thus,
    \[\log F(qx)-\log F(x)=\log(1-x)=-\sum_{k=1}^\infty\dfrac{x^k}{k}.\]
    Let $a_k$ be the coefficient of the $x^k$-term in the power series $\log F(x)$. Then
    \[(q^k-1)a_k=-\dfrac{1}{k}.\]
    So the right-hand side of the equality we are proving is
    \[\dfrac{1}{r}\prod_{j=1}^{r-1}(1-q^j)=d_ra_r.\]
    If we calculate $a_r$ directly, then we get
    \[a_r=\sum_{r_1+\dots+r_k=r}\frac{(-1)^{k-1}}{k}\frac{1}{d_{r_1}\cdots d_{r_k}}.\]
    So
    \begin{align*}
        \RHS=d_ra_r&=\sum_{r_1+\dots+r_k=r}\frac{(-1)^{k-1}}{k}\frac{d_r}{d_{r_1}\cdots d_{r_k}}\\
        &=\sum_{r_1+\dots+r_k=r}\dfrac{(-1)^{k-1}}{k}\frac{\displaystyle{\prod_{j=1}^{r}(1-q^j)}}{\displaystyle{\prod_{j=1}^{r_1}(1-q^j)\dots\prod_{j=1}^{r_k}(1-q^j)}}\\
        &=\sum_{r_1+\dots+r_k=r}\dfrac{(-1)^{k-1}}{k}\frac{\displaystyle{\prod_{j=1}^{r}(q^j-1)}}{\displaystyle{\prod_{j=1}^{r_1}(q^j-1)\dots\prod_{j=1}^{r_k}(q^j-1)}}=\LHS.
    \end{align*}
\end{proof}

Since the coefficient of $u_{S_i^{\oplus r}[b]}$ in $\epsilon^{\sigma_*}_{\phi_*(S_i)+b,(-1)^b re_i}$   is not zero, we have
\[\mathfrak{L}_{\sigma_*}=\Lie_{\mathbb{Q}}\linp{u_{S_i^{\oplus r}[b]}\middle|i\in Q_0,r\in\mathbb{Z}_+,b\in\mathbb{Z}}\subseteq\mathcal{H}(\mathcal{D}_Q).\]

Combining the preceding lemmas with Theorem \ref{finite length Joyce transport} gives the following results:
\begin{theorem}[Joyce transport for Dynkin cases]
    For any $\sigma\in\Stab(\mathcal{D}_Q),\gamma\in\mathbb{Z}^n,\phi\in\mathbb{R}$, the Joyce element $\epsilon^\sigma_{\phi,\gamma}$ can be expressed as a Lie polynomial in $u_{S_i^{\oplus r}[b]}$ with coefficients in $\mathbb{Q}$, where $i\in Q_0,r\in\mathbb{Z}_+,b\in\mathbb{Z}$. As a consequence, the Lie subalgebra at $\sigma$ is
    \[\mathfrak{L}_\sigma=\mathfrak{L}_{\sigma_*}=\Lie_{\mathbb{Q}}\linp{u_{S_i^{\oplus r}[b]}\middle|i\in Q_0,r\in\mathbb{Z}_+,b\in\mathbb{Z}}\subseteq\mathcal{H}(\mathcal{D}_Q).\]
\end{theorem}
 
\section*{Acknowledgments}
This work was supported by NSF of China (No. 12371036).

\section*{Declaration of AI use}
The authors used ChatGPT to assist in improving some of the mathematical arguments and the English expression. All mathematical results, arguments, and proofs were carefully developed and verified by the authors.

\end{document}